\documentclass[graybox]{svmult}

\usepackage{mathptmx}       
\usepackage{helvet}         
\usepackage{courier}        
\usepackage{type1cm}        
\usepackage{makeidx}         
\usepackage{graphicx}        
\usepackage{multicol}        
\usepackage[bottom]{footmisc}
\usepackage{amsmath,amsfonts,amssymb,pifont}

\newcommand{\be}[1]{\begin{equation}\label{#1}}
\newcommand{\ee}{\end{equation}}
\newcommand{\ba}[1]{\begin{eqnarray}\label{#1}}
\newcommand{\ea}{\end{eqnarray}}
\newcommand{\rf}[1]{(\ref{#1})}
\newcommand{\nn}{\nonumber}

\makeindex             

\begin{document}

\title*{From rotating fluid masses and Ziegler's paradox to
Pontryagin- and Krein spaces and bifurcation theory
}
\titlerunning{Rotating fluid masses, Ziegler's paradox, Krein spaces and bifurcation theory}
\author{Oleg N. Kirillov and Ferdinand Verhulst}
\institute{Oleg N. Kirillov \at Northumbria University, NE1 8ST Newcastle upon Tyne, UK, \email{oleg.kirillov@northumbria.ac.uk}
\and Ferdinand Verhulst \at Mathematisch Instituut, PO Box 80.010, 3508TA Utrecht, Netherlands \email{F.Verhulst@uu.nl}}
%
%
\maketitle

\abstract*{Three classical systems, the Kelvin gyrostat, the Maclaurin spheroids, and the Ziegler pendulum have directly inspired development of the theory of Pontryagin and Krein spaces with indefinite metric and singularity theory as independent mathematical topics, not to mention stability theory and nonlinear dynamics. From industrial applications in shipbuilding, turbomachinery, and artillery to fundamental problems of astrophysics, such as asteroseismology and gravitational radiation --- that is the range of phenomena involving the Krein collision of eigenvalues, dissipation-induced instabilities, and spectral and geometric singularities on the neutral stability surfaces, such as the famous Whitney's umbrella.}

\abstract{Three classical systems, the Kelvin gyrostat, the Maclaurin spheroids, and the Ziegler pendulum have directly inspired development of the theory of Pontryagin and Krein spaces with indefinite metric and singularity theory as independent mathematical topics, not to mention stability theory and nonlinear dynamics. From industrial applications in shipbuilding, turbomachinery, and artillery to fundamental problems of astrophysics, such as asteroseismology and gravitational radiation --- that is the range of phenomena involving the Krein collision of eigenvalues, dissipation-induced instabilities, and spectral and geometric singularities on the neutral stability surfaces, such as the famous Whitney's umbrella.}

\section{Historical background}
\label{sec:1}
The purpose of this paper is to show how  a curious phenomenon observed in the natural sciences, destabilization by
dissipation, was solved by mathematical analysis. After the completion of the analysis, the eigenvalue
calculus of matrices which it involved was (together with other applications) an inspiration
for the mathematical  theory of structural stability of matrices.  But the story is even more intriguing.
Later it was shown that the bifurcation picture in parameter space is related to a seemingly pure mathematical object in singularity theory, Whitney's umbrella.

In many problems in physics and engineering, scientists found what looked like a counter-intuitive
phenomenon: certain systems that without dissipation show stable behaviour became unstable when any form of dissipation was
introduced. This can be viscosity in fluids and solids, magnetic diffusivity, or losses due to radiation of waves of different nature,
including recently detected gravitational waves \cite{Abb17,A2003}, to name a few. For instance, dissipation-induced modulation instabilities \cite{BD2007} are widely discussed in the context of modern nonlinear optics \cite{PTS2018}.

The destabilization by dissipation is especially sophisticated when several dissipative mechanisms are acting simultaneously.  In  this  case,  ``no  simple  rule  for  the  effect  of  introducing  small viscosity  or  diffusivity  on  flows  that  are  neutral  in  their  absence  appears  to  hold'' \cite{TSL2013} and ``the ideal limit with zero dissipation coefficients has essentially nothing to do with the case of small but finite dissipation coefficients'' \cite{M1993}.

In hydrodynamics,  a  classical  example  is  given  by  secular  instability  of  the  Maclaurin  spheroids due to both fluid viscosity (Kelvin and Tait, 1879) and gravitational radiation reaction (Chandrasekhar, 1970), where the critical eccentricity of the meridional section of the spheroid depends on the ratio of the two dissipative mechanisms and  reaches  its  maximum,  corresponding  to  the  onset  of  dynamical  instability  in  the  ideal system, when  this  ratio  equals  1  \cite{LD1977}. In meteorology this phenomenon is known as the `Holop\"ainen instability mechanism' (Holop\"ainen, 1961) for a
baroclinic flow when waves that are linearly stable in the absence of Ekman friction become
dissipatively destabilized in its presence, with the result that the location of the curve of marginal
stability is displaced by an order one distance in the parameter space, even if the Ekman number
is infinitesimally small \cite{Ho61,WE2012}. For a baroclinic circular vortex with thermal and viscous diffusivities this phenomenon was studied by McIntyre in 1970 \cite{MI1970}.

In  rotor dynamics,  the  generic  character of  the  discontinuity  of  the  instability  threshold  in  the  zero  dissipation  limit  was  noticed by  Smith already in 1933 \cite{S33}.
In mechanical engineering such a phenomenon is called Ziegler's paradox, it was found in the analysis of a double pendulum with a nonconservative positional force with and without damping in 1952 \cite{Z52,Z53}.
The importance of solving the Ziegler paradox for mechanics was emphasized by Bolotin \cite{Bo63}: ``The greatest theoretical interest is evidently centered in the unique effect of damping in the presence of pseudo-gyroscopic forces, and
in particular, in the differences in the results for systems with slight damping which then becomes zero and systems in which damping
is absent from the start.'' Encouraging further research of the destabilization paradox, Bolotin was not aware that the crucial ideas for its explanation were formulated as early as 1956.

Ziegler's paradox was solved in 1956 by an expert in classical geometry and mechanics, Oene Bottema \cite{B56}. He formulated
the problem of the stability of an equilibrium in two degrees-of-freedom (4 dimensions), allowing for gyroscopic
and nonconservative positional forces. The solution by Bottema in the form of concrete analysis was hardly noticed at that time. Google
Scholar gives no citation of Bottema's paper in the first 30 years after 1956.

As mentioned above,
a new twist to the treatment of the problem came from identifying geometric considerations independently introduced
by Whitney in singularity theory with the bifurcation analysis of Bottema. Interestingly Whitney's ``umbrella singularity''
predated Bottema's analysis, it gives the right geometric picture but it ignores the stability questions of the dynamical
context which is the essential question of Ziegler's paradox. Later, in 1971-72,
V.I. Arnold showed that the umbrella singularity is generic in parameter families of real matrices. This result links
to stability by linearisation of the vector field near equilibrium.

The phenomenon of dissipation-induced instability and the Ziegler paradox raised important questions in mechanics and mathematics.
For instance, what is the connection between the conservative and deterministic Hamiltonian systems and real systems involving both dissipation of energy and stochastic effects \cite{La03}? Also it added a new bifurcation in the analysis of dynamical systems. An important consequence of the results is that in a large number of problem-fields one can now predict and characterise precisely this type of instability \cite{HR95, KV10, K2013dg, KI2017, KM07, La03, LFA2016, SJ07}.

\subsection{Stability of Kelvin's gyrostat and spinning artillery shells filled with liquid}

      \begin{figure}
    \begin{center}
    \includegraphics[angle=0, width=0.4\textwidth]{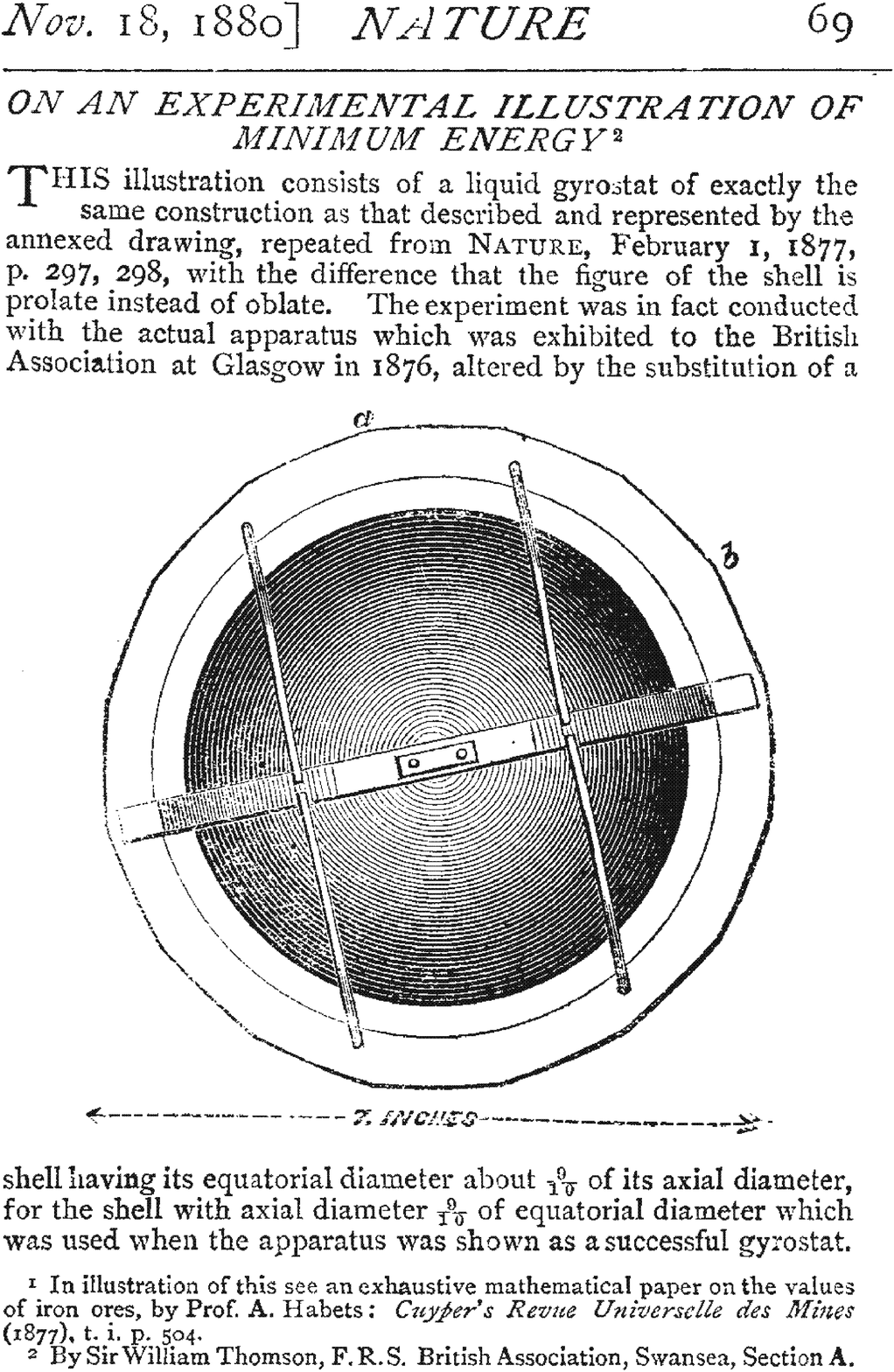}
    \includegraphics[angle=0, width=0.4\textwidth]{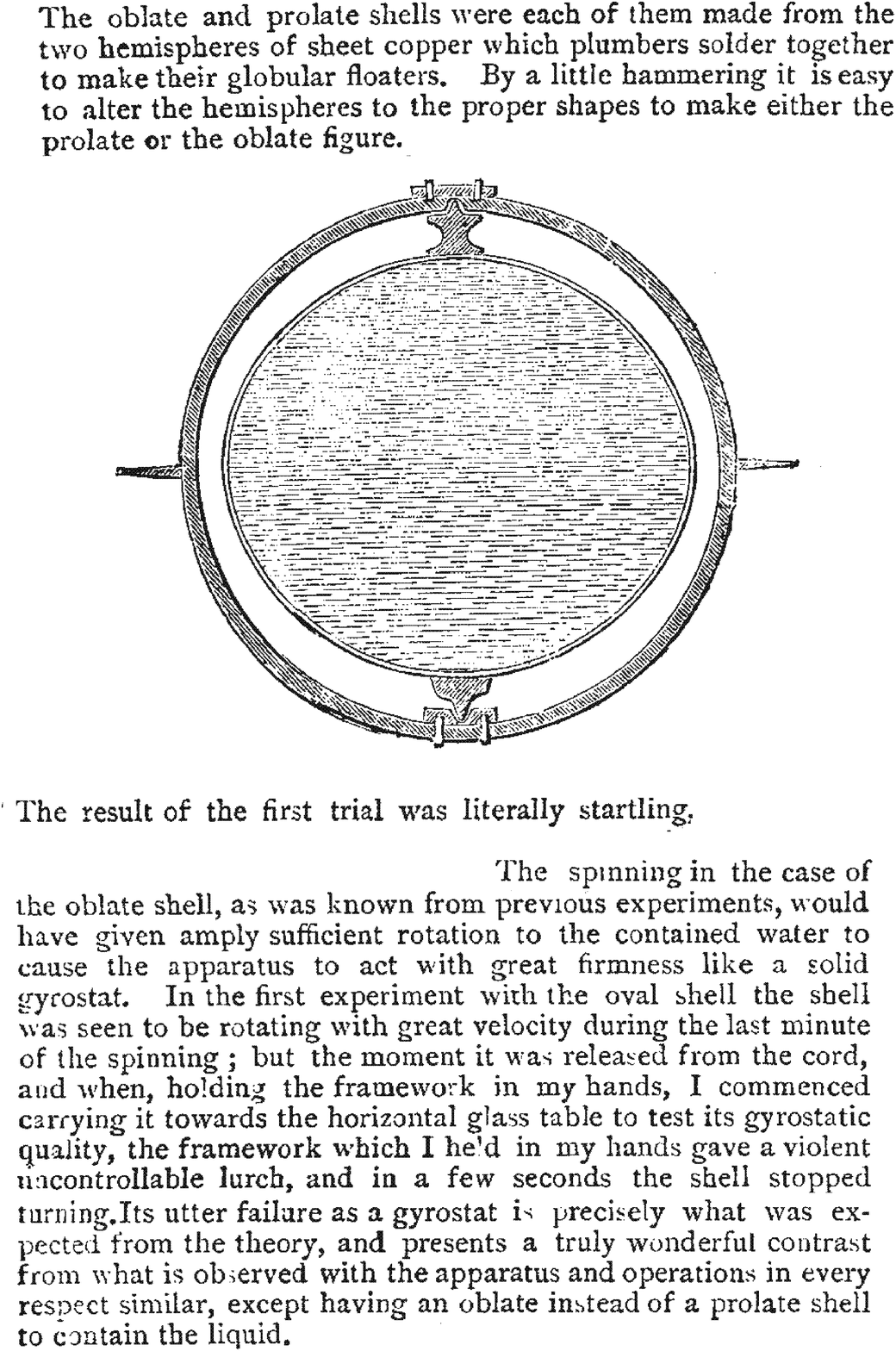}
    \end{center}
    \caption{The Kelvin gyrostat \cite{WT1880}.}
    \label{kelgyro}
    \end{figure}

Already von Laue (1905) \cite{VL1905}, Lamb (1908) \cite{Lamb1908} and Heisenberg (1924) \cite{H1924} realized that dissipation easily destabilizes waves and modes of negative energy of an ideal system supported by rotating and translating continua \cite{S1958,S1960}. Williamson (1936) proposed normal forms for Hamiltonian systems allowing sorting stable modes of negative and positive energy according to their symplectic sign \cite{MK86,W1936,W1937}. However, further generalization --- the theory of spaces with indefinite metric, or Krein \cite{K1983} and Pontryagin \cite{P44} spaces, --- was directly inspired by the problem of stability of the Kelvin gyrostat \cite{G1880,WT1880}, having both astrophysical and industrial (and even military) applications.

Kelvin \cite{WT1880} experimentally demonstrated in 1880
that a thin-walled and slightly oblate spheroid completely filled with liquid
remains stable if rotated fast enough about a fixed point, which does not happen if the
spheroid is slightly prolate, Figure~\ref{kelgyro}. In the same year this observation was confirmed theoretically
by Greenhill \cite{G1880}, who found that rotation around the center of gravity of the top in
the form of a weightless ellipsoidal shell completely filled with an ideal and incompressible
fluid, is unstable when $a < c < 3a$, where $c$ is the length of the semiaxis of
the ellipsoid along the axis of rotation and the lengths of the two other semiaxes are
equal to $a$ \cite{G1880}.

Quite similarly, bullets and projectiles fired from the rifled weapons can relatively easily be stabilized by rotation,
if they are solid inside. In contrast, the shells, containing a liquid substance inside, have a tendency to turn over
despite seemingly revolved fast enough to be gyroscopically stabilized.
Motivated by such artillery applications, in 1942 Sobolev,
then director of the Steklov Mathematical Institute in Moscow, considered stability of a rotating heavy
top with a cavity entirely filled with an ideal incompressible fluid \cite{S2006}---a problem
that is directly connected to the classical XIXth century models of astronomical bodies with a crust surrounding a molten
core \cite{S1959}.

For simplicity, the solid shell of the top and the domain $V$ occupied by the cavity inside it, can be assumed to
have a shape of a solid of revolution. They have a common symmetry axis where the fixed point of the top is located.
The velocity profile of the stationary unperturbed motion of the fluid is that of a solid body rotating with the same
angular velocity $\Omega$ as the shell around the symmetry axis.

Following Sobolev, we denote by $M_1$ the mass of the shell, $M_2$ the mass of the fluid,
$\rho$ and $p$ the density and the pressure of the fluid, $g$ the gravity acceleration, and $l_1$ and $l_2$ the distances from the fixed point to the centers of mass of the shell and the fluid, respectively. The moments of inertia of the shell and the `frozen' fluid  with respect to the symmetry axis are $C_1$ and $C_2$, respectively;
$A_1$ ($A_2$) stands for the moment of inertia of the shell (fluid) with respect to any axis that is orthogonal to the symmetry axis and passes through the fixed point. Let, additionally,
\be{sob1}
L=C_1+C_2-A_1-A_2-\frac{K}{\Omega^2},\quad K=g(l_1M_1+l_2M_2).
\ee
The solenoidal $(\textrm{div}\, \textbf{v}=0)$ velocity field $\textbf{v}$ of the fluid is assumed to satisfy the no-flow condition on the boundary of the cavity: $\left. \textbf{ v}_n \right|_{\partial V}=0$.

Stability of the stationary rotation of the top around its vertically oriented symmetry axis
is determined by the system of linear equations derived by Sobolev in the frame $(x,y,z)$
that has its origin at the fixed point of the top and rotates with respect to an inertial frame around the
vertical $z$-axis with the angular velocity of the unperturbed top, $\Omega$. If the real and imaginary part of the
complex number $Z$ describe the deviation of the unit vector of the symmetry axis of the top in the coordinates $x$, $y$, and $z$,
then these equations are, see e.g. \cite{S2006,Y97}:
\ba{sob}
\frac{d Z}{d t}&=&i \Omega W,\nn\\
(A_1{+}\rho \kappa^2)\frac {d W}{dt}&=&i\Omega L Z+i\Omega(C_1{-}2A_1 {+}\rho E) W\nn+ i\rho
\int_V\left(v_x \frac{\partial \chi}{\partial y}-v_y \frac{\partial \chi}{\partial x} \right)dV,\nn\\
\partial_t v_x&=&2\Omega v_y-{\rho}^{-1}\partial_x p+2i\Omega^2 W\partial_y \overline \chi,\nn\\
\partial_t v_y&=&-2\Omega v_x-{\rho}^{-1}\partial_y p-2i\Omega^2 W\partial_x \overline \chi,\nn\\
\partial_t v_z&=&-{\rho}^{-1}\partial_z p,
\ea
where $2\kappa^2=\int_V |\nabla \chi|^2 dV$, $E=i\int_V\left(\partial_x \overline \chi \partial_y \chi-\partial_y \overline \chi \partial_x \chi  \right) dV$, and the function $\chi$ is determined
by the conditions
\be{sob1}
\nabla^2 \chi=0,\quad \left. \partial_n \chi \right|_{\partial V}=z(\cos n x+i \cos n y)-(x+iy)\cos n z,
\ee
with $n$ the absolute value of a vector $\textbf{n}$, normal to the boundary of the cavity.

Sobolev realized that some qualitative conclusions on the stability of the top can be drawn with the use of the
bilinear form
\be{qform}
Q(R_1,R_2)=L\Omega Z_1\overline Z_2+(A_1+\rho \kappa^2)W_1\overline W_2+\frac{\rho}{2\Omega^2}\int_V\overline{\textbf{v}}_2^T\textbf{v}_1dV
\ee
on the elements $R_1$ and $R_2$ of the space $\{R\}=\{Z,W,\textbf{v}\}$. The linear operator $B$ defined by Eqs.~\rf{sob}
that can be written as $\frac{ d R}{d t}=iBR$ has all its eigenvalues real when $L>0$, which yields Lyapunov stability of
the top. The number of pairs of complex-conjugate eigenvalues of $B$ (counting multiplicities) does not exceed the number
of negative squares of the quadratic form $Q(R,R)$, which can be equal only to one when $L<0$. Hence, for $L<0$ an unstable solution
$R=e^{i\lambda_0 t}R_0$ can exist with $\textrm{ Im}\lambda_0<0$; all real eigenvalues are simple except for maybe one
\cite{S2006,Y97,Kopachevskii2001}.

In the particular case when the cavity is an ellipsoid of rotation with the semi-axes $a$, $a$, and $c$,
the space of the velocity fields of the fluid can be decomposed into a direct sum of subspaces, one of which is finite-dimensional. Only the movements from this subspace interact with the movements of the rigid shell,
which yields a finite-dimensional system of ordinary differential equations that describes coupling between the shell and the fluid.

Calculating the moments of inertia of the fluid in the ellipsoidal container
$$
C_2=\frac{8\pi\rho}{15} a^4c,\quad A_2=l_2^2M_2+\frac{4\pi\rho}{15}  a^2 c (a^2+c^2),
$$
denoting $m=\frac{c^2-a^2}{c^2+a^2}$, and assuming the field $\textbf{ v}=(v_x,v_y,v_z)^T$ in the form
$$
v_x=(z-l_2)a^2m\xi,\quad v_y=-i(z-l_2)a^2m\xi,\quad v_z=-(x-iy)c^2m\xi,
$$
one can eliminate the pressure in Eqs.~\rf{sob} and obtain the  reduced model
\be{meqs}
\frac{d\textbf{x}}{dt}=i\Omega \textbf{A}^{-1}\textbf{C}\textbf{x}=i\Omega\textbf{B}\textbf{x},
\ee
where $\textbf{x}=(Z,W,\xi)^T \in \mathbb{C}^3$ and
\ba{sobell}
\textbf{A}&=&\left(
  \begin{array}{ccc}
    1 & 0 & 0 \\
    0 & A_1{+}l_2^2M_2{+}\frac{4\pi \rho}{15}a^2c\frac{(c^2-a^2)^2}{c^2+a^2} & 0 \\
    0 & 0 & c^2+a^2 \\
  \end{array}
\right),\nn\\
\textbf{C}&=&\left(
  \begin{array}{ccc}
    0 & 1 & 0 \\
     L & C_1{-}2A_1{-}2l_2^2M_2{-}\frac{8\pi \rho}{15}a^2c^3m^2& -\frac{8\pi \rho}{15}a^4c^3m^2   \\
    0 & -2  & -2  a^2 \\
  \end{array}
\right).
\ea

The matrix $\textbf{B}\ne \textbf{B}^T$ in Eq.~\rf{meqs} after multiplication by a symmetric matrix
\be{matg}
\textbf{G}=\left(
          \begin{array}{ccc}
            L & 0 & 0 \\
            0 & A_1{+}l_2^2M_2{+}\frac{4\pi \rho}{15}a^2c\frac{(c^2-a^2)^2}{c^2+a^2} & 0 \\
            0 & 0 & \frac{4\pi \rho}{15}a^4c^3\frac{(c^2-a^2)^2}{c^2+a^2} \\
          \end{array}
        \right)
\ee
yields a Hermitian matrix $\textbf{G}\textbf{B}=\overline{(\textbf{G}\textbf{B})}^T$, i.e. $\textbf{B}$ is a self-adjoint operator in the
space $\mathbb{C}^3$ endowed with the metric
\be{indefp}
[\textbf{u}, \textbf{u}]:=(\textbf{G}\textbf{u},\textbf{u})=\overline{\textbf{u}}^T\textbf{G}\textbf{u},\quad \textbf{u} \in \mathbb{C}^3,
\ee
which is \emph{definite} when $L>0$ and \index{indefinite metric} \emph{indefinite} with one negative square when $L<0$. If $\lambda$ is
an eigenvalue of the matrix $\textbf{B}$, i.e. $\textbf{B}\textbf{u}=\lambda \textbf{u}$,
then $\overline{\textbf{u}}^T\textbf{G}\textbf{B}\textbf{u}=\lambda\overline{\textbf{u}}^T\textbf{G}\textbf{u}$. On the other hand,
$\overline{\textbf{u}}^T (\overline{\textbf{G}\textbf{B}})^T\textbf{u}=\overline \lambda\,
\overline{\overline{\textbf{u}}^T\textbf{G}\textbf{u}}=\overline \lambda\, {\overline{\textbf{u}}^T\textbf{G}\textbf{u}}$. Hence,
$$
(\lambda-\overline\lambda)\overline{\textbf{u}}^T\textbf{G}\textbf{u}=0,
$$
implying $\overline{\textbf{u}}^T\textbf{G}\textbf{u}=0$ on the eigenvector $\textbf{u}$ of the complex $\lambda\ne \overline \lambda$.
For real eigenvalues $\lambda= \overline \lambda$ and
$\overline{\textbf{u}}^T\textbf{G}\textbf{u}\ne 0$. The sign of the quantity $\overline{\textbf{u}}^T\textbf{G}\textbf{u}$ (or Krein sign) can be
different for different real eigenvalues.

    \begin{figure}[htp]
    \begin{center}
    \includegraphics[angle=0, width=0.95 \textwidth]{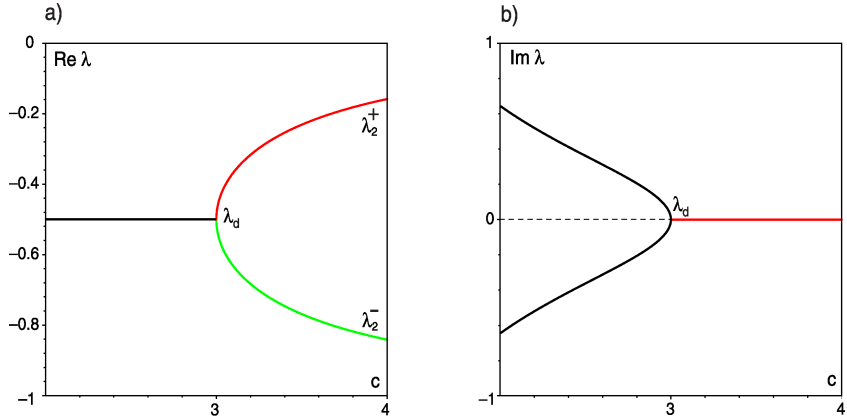}
    \end{center}
    \caption{(a) Simple real eigenvalues \rf{com} of the Sobolev's top in the Greenhill's case for $a=1$ with (red)
    $\overline{\textbf{u}}^T\textbf{G}\textbf{u}>0$ and (green) $\overline{\textbf{u}}^T\textbf{G}\textbf{u}<0$. (b)
    At simple complex-conjugate eigenvalues (black) and at the double real eigenvalue $\lambda_d$ we have
    $\overline{\textbf{u}}^T\textbf{G}\textbf{u}=0$.}
    \label{fig03_1}
    \end{figure}

For example, when the ellipsoidal shell is massless and the supporting point is at the center of mass of the system, then
$A_1=0$, $C_1=0$, $M_1=0$, $l_2=0$. The matrix $\textbf{B}$ has thus one real eigenvalue
$(\lambda_1^+=-1,~ {\overline{\textbf{u}}_{1}^+}^T\textbf{G}\textbf{u}_{1}^+>0)$ and the pair of eigenvalues
\be{com}
\lambda_2^{\pm}=-\frac{1}{2}\pm\frac{1}{2}\sqrt{1+\frac{32\pi\rho}{15}\frac{ca^4}{L}},\quad
L=\frac{4\pi \rho}{15}a^2c(a^2-c^2),
\ee
which are real if $L>0$ and can be complex if $L<0$. The latter condition together with the requirement that the radicand in Eq.~\rf{com}
is negative, reproduces the Greenhill's instability zone: $a<c<3a$ \cite{G1880}.
With the change of $c$, the  real eigenvalue $\lambda_2^{+}$ with ${\overline{\textbf{u}}_{2}^+}^T\textbf{G}\textbf{u}_{2}^+>0$ collides
at $c=3a$ with the real eigenvalue $\lambda_2^{-}$ with ${\overline{\textbf{u}}_{2}^-}^T\textbf{G}\textbf{u}_{2}^-<0$ into a real
double defective eigenvalue $\lambda_d$ with the algebraic multiplicity two and geometric multiplicity one. This Krein collision is illustrated in Figure~\ref{fig03_1}.
Note that ${\overline{\textbf{u}}_{d}}^T\textbf{G}\textbf{u}_{d}=0$, where $\textbf{u}_d$ is the eigenvector at $\lambda_d$.

Therefore, in the case of the ellipsoidal shapes of the shell and the cavity, the Hilbert space $\{R\}=\{Z,W,\textbf{v}\}$ of the
Sobolev's problem endowed with the indefinite metric $(L<0)$ decomposes into the three-dimensional space of the reduced model
\rf{meqs}, where the self-adjoint operator $B$ can have complex eigenvalues and real defective eigenvalues, and a complementary
infinite-dimensional space, which is free of these complications. The very idea that the signature of the indefinite metric can
serve for counting unstable eigenvalues of an operator that is self-adjoint in a functional space equipped with such a metric,
turned out to be a concept of a rather universal character possessing powerful generalizations
that were initiated by Pontryagin in 1944 \cite{P44} and developed into a general theory of indefinite inner product spaces or Krein spaces \cite{Bognar1974,KM2014,K1983}. Relation of the Krein sign to the sign of energy or action has made it a popular tool for predicting instabilities in physics \cite{N1976,Z2016}.

\subsection{Secular instability of the Maclaurin spheroids by viscous and radiative losses}

It is hard to find a physical application that would stimulate development of mathematics to such an extent as the problem
of stability of equilibria of rotating and self-gravitating masses of fluids. Rooted in the Newton and Cassini thoughts on the actual shape of the Earth, the rigorous analysis of this question attracted the best minds of the XVIII-th and XIX-th centuries, from Maclaurin to Riemann, Poincar\'e and Lyapunov.
In fact, modern nonlinear dynamics \cite{AA88,IA92,Kuz04} and Lyapunov stability theory \cite{L1992} are by-products of the efforts invested in solution of this question of the astrophysical fluid dynamics \cite{BKM2009}, which experiences a revival nowadays \cite{A2003} inspired by the recent detection of gravitational waves \cite{Abb17}.

We recall that in 1742 Maclaurin established that an oblate spheroid
$$
\frac{x^2}{a_1^2}+\frac{y^2}{a_2^2}+\frac{z^2}{a_3^2}=1, \quad a_3<a_2=a_1
$$
is a shape of relative equilibrium of a self-gravitating mass of
inviscid fluid in a solid-body rotation about the $z$-axis, provided that the rate of rotation, $\Omega$, is
related to the eccentricity  $e=\sqrt{1-\frac{a_3^2}{a_1^2}}$ through the formula \cite{M1742}
\be{maclaurin}
\Omega^2(e)=2e^{-3}(3-2e^2)\sin^{-1}(e)\sqrt{1-e^2}-6e^{-2}(1-e^2).
\ee

A century later, Jacobi (1834) has discovered less symmetric shapes of relative equilibria in this problem that are tri-axial ellipsoids
$$
\frac{x^2}{a_1^2}+\frac{y^2}{a_2^2}+\frac{z^2}{a_3^2}=1, \quad a_3<a_2<a_1.
$$
Later on a student of Jacobi, Meyer (1842) \cite{M1842}, and then Liouville (1851) \cite{L1851} have shown that  the family of Jacobi's
ellipsoids has one member in common with the family of Maclaurin's spheroids at $e\approx 0.8127$.
The equilibrium with the Meyer-Liouville eccentricity is neutrally stable.

In 1860 Riemann \cite{R1860} established neutral stability of inviscid Maclaurin's spheroids on the interval of eccentricities $(0<e<0.952..)$.
At the Riemann point with the critical eccentricity $e\approx0.9529$ the Hamilton-Hopf bifurcation sets in and causes
dynamical instability with respect to ellipsoidal perturbations beyond this point \cite{DS1979,S1980-2,S1980-3}.

A century later Chandrasekhar \cite{C1969} used a virial theorem to reduce the problem to a finite-dimensional system,
which stability is governed by the eigenvalues of the matrix polynomial
\be{mpi}
{\mathbf L}_i(\lambda)=\lambda^2\left(
           \begin{array}{cc}
             1 & 0 \\
             0 & 1 \\
           \end{array}
         \right)+\lambda\left(
                          \begin{array}{cc}
                            0 & -4\Omega \\
                            \Omega & 0 \\
                          \end{array}
                        \right)+
                        \left(
                          \begin{array}{cc}
                            4b-2\Omega^2 & 0 \\
                            0 & 4b-2\Omega^2 \\
                          \end{array}
                        \right),
\ee
where $\Omega(e)$ is given by the Maclaurin law \rf{maclaurin} and  $b(e)$ is as follows
\be{macb}
b = \frac{\sqrt{1-e^2}}{4e^5}\left\{e(3-2e^2)\sqrt{1-e^2}+(4e^2-3)\sin^{-1}(e)\right\}.
\ee
The eigenvalues of the matrix polynomial \rf{mpi} are
\be{eigmp}
\lambda=\pm \left(i \Omega\pm i \sqrt{4b-\Omega^2}\right).
\ee

Requiring $\lambda=0$ we can determine the critical Meyer-Liouville eccentricity by solving with respect to $e$ the equation \cite{C1969}
$$
4b(e)=2\Omega^2(e).
$$
The critical eccentricity at the Riemann point follows from requiring the radicand in \rf{eigmp} to vanish:
$$
4b(e)=\Omega^2(e),
$$
which is equivalent to the equation
$$
e=\sin\left(\frac{e(3+4e^2)\sqrt{1-e^2}}{3+2e^2-4e^4}\right)
$$
that has a root $e\approx0.9529$.

\begin{figure}
\begin{center}
\includegraphics[angle=0, width=0.47\textwidth]{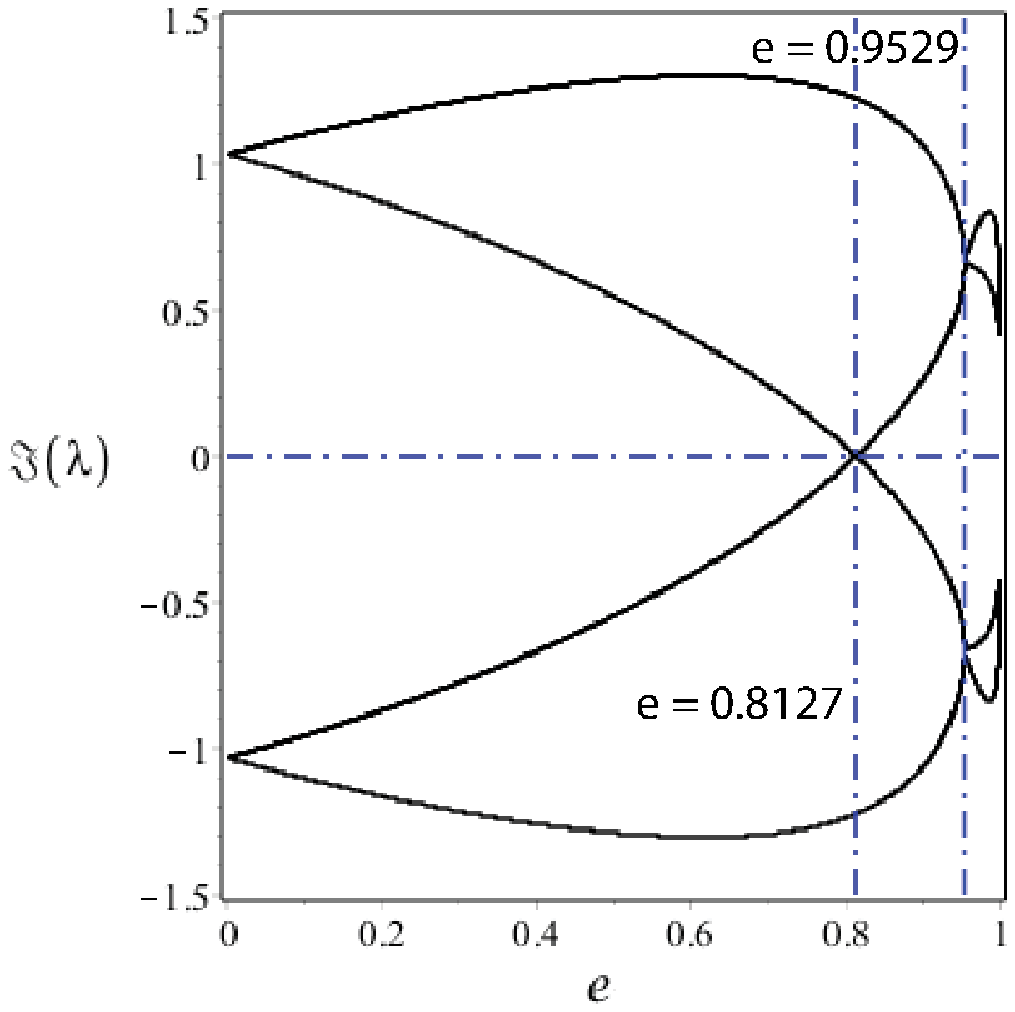}
\includegraphics[angle=0, width=0.47\textwidth]{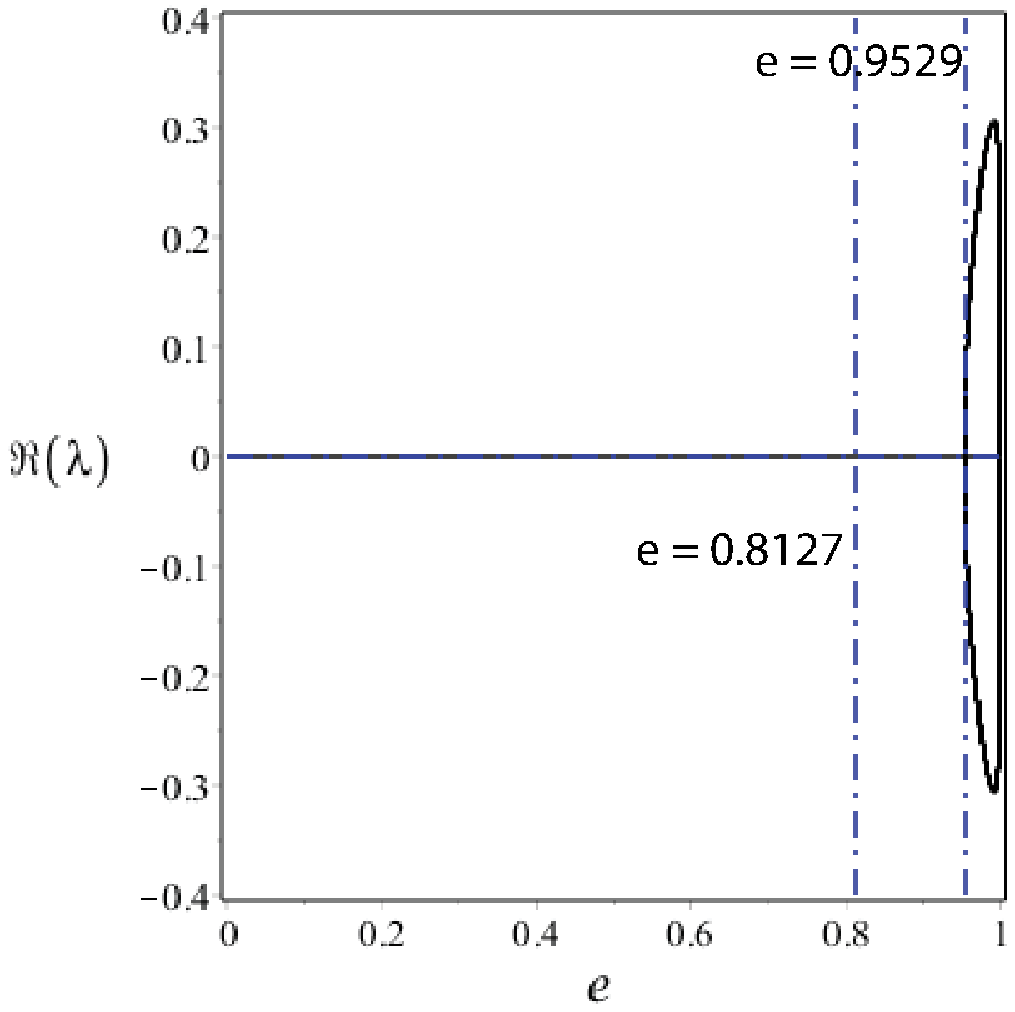}
\end{center}
\caption{(Left) Frequencies and (right) growth rates of the eigenvalues of the inviscid eigenvalue problem $\textbf{L}_i(\lambda)\textbf{u}=0$ demonstrating
the Hamilton-Hopf bifurcation at the Riemann critical value of the eccentricity, $e \approx 0.9529$ and neutral stability at the Meyer-Liouville point,
$e\approx0.8127$.}
\label{figcond1}
\end{figure}

Remarkably, when
\be{interval}
\Omega^2(e)<4b(e)<2\Omega^2(e)
\ee
both eigenvalues of the stiffness matrix
$$
\left(
                          \begin{array}{cc}
                            4b-2\Omega^2 & 0 \\
                            0 & 4b-2\Omega^2 \\
                          \end{array}
                        \right)
$$
are negative. The number of negative eigenvalues of the matrix of potential forces is known as the Poincar\'e instability degree.
The Poincar\'e instability degree of the equilibria with the eccentricities \rf{interval} is even and equal to 2.
Hence, the interval \rf{interval} corresponding to $0.812..<e<0.952..$, which is stable according to Riemann,
is, in fact, the interval of gyroscopic stabilization \cite{L2013} of the Maclaurin spheroids, Figure~\ref{figcond1}.

In 1879 Kelvin and Tait \cite{TT} realized that viscosity of the fluid can destroy the gyroscopic stabilization of the Maclaurin spheroids: ``If there be any viscosity, however slight,
in the liquid, the equilibrium [beyond $e\approx 0.8127$] in any case of energy either a minimax or a maximum cannot be secularly stable''.

The prediction made by Kelvin and Tait \cite{TT} has been rigorously verified only in the XX-th century by Roberts and Stewartson \cite{RS1963}. Using the virial approach
by Chandrasekhar, the linear stability problem can be reduced to the study of eigenvalues of the matrix polynomial
\be{mprs}
{\mathbf L}_{v}(\lambda)=\lambda^2\left(
           \begin{array}{cc}
             1 & 0 \\
             0 & 1 \\
           \end{array}
         \right)+\lambda\left(
                          \begin{array}{cc}
                            10\mu & -4\Omega \\
                            \Omega & 10\mu \\
                          \end{array}
                        \right)+
                        \left(
                          \begin{array}{cc}
                            4b-2\Omega^2 & 0 \\
                            0 & 4b-2\Omega^2 \\
                          \end{array}
                        \right),
\ee
where $\mu=\frac{\nu}{a_1^2\sqrt{\pi G \rho}}$, $\nu$ is the viscosity of the fluid, $G$ is the universal gravitation constant, and $\rho$ the density of the fluid \cite{C1969}. The $\lambda$ and $\Omega$ are measured in units of $\sqrt{\pi G \rho}$. The operator ${\mathbf L}_{v}(\lambda)$ differs
from the operator of the ideal system, ${\mathbf L}_{i}(\lambda)$, by the matrix of dissipative forces $10\lambda\mu \textbf{I}$, where
$\textbf{I}$ is the $2\times 2$ unit matrix.

\begin{figure}
\begin{center}
\includegraphics[angle=0, width=0.47\textwidth]{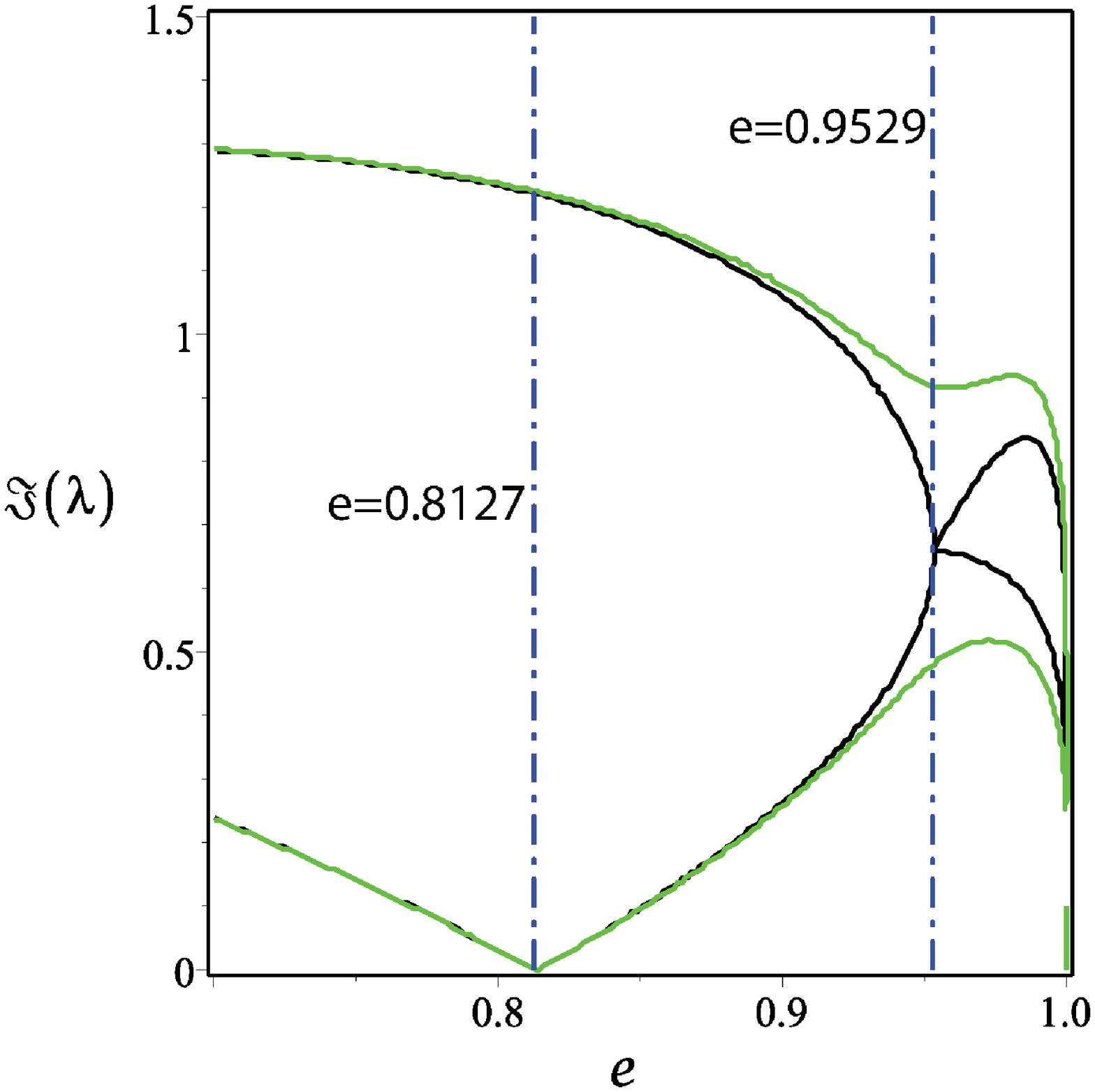}
\includegraphics[angle=0, width=0.47\textwidth]{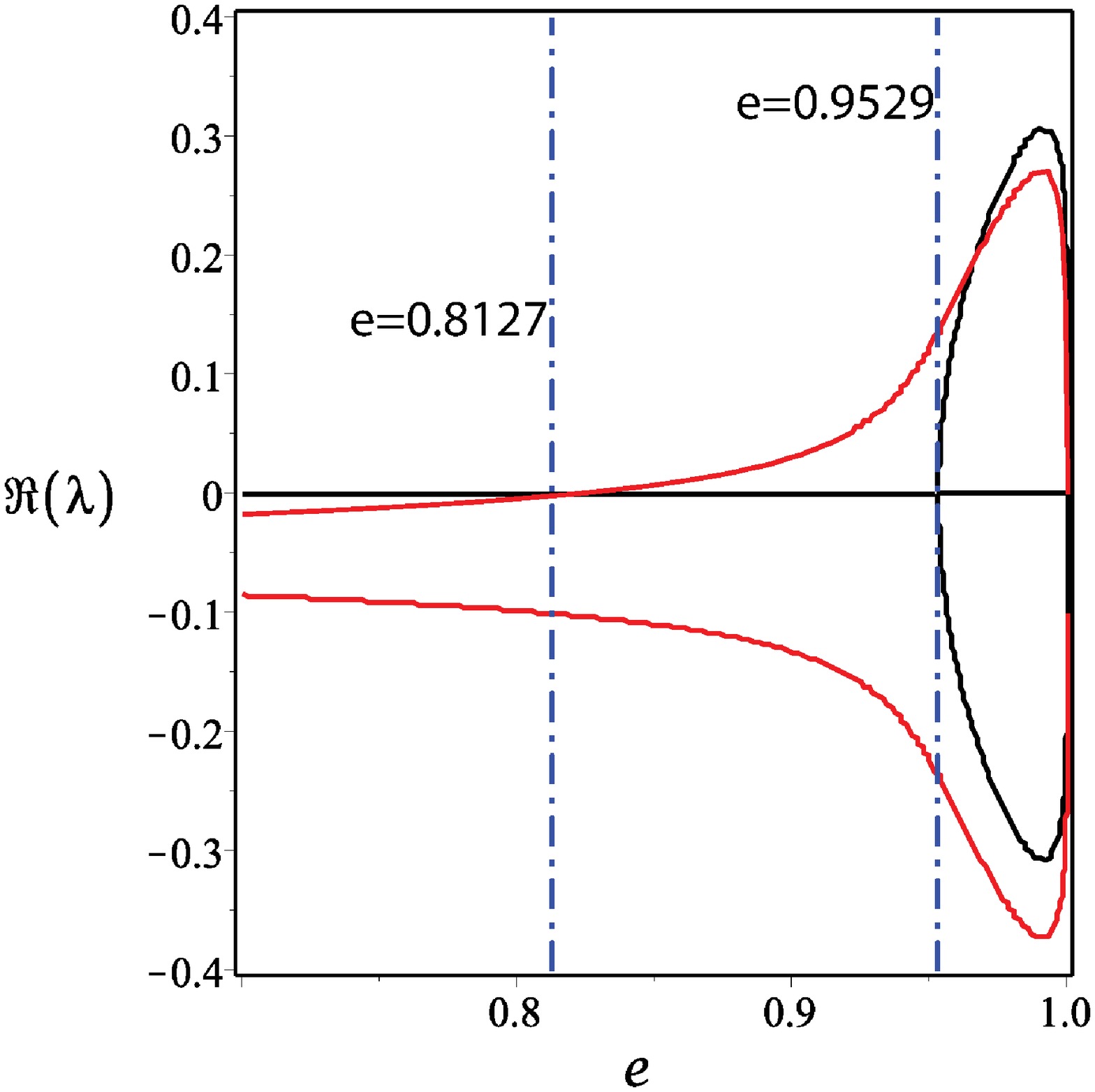}
\end{center}
\caption{(Left) Frequencies and (right) growth rates of the (black lines) inviscid Maclaurin spheroids and (green and red lines) viscous
ones with $\mu=\frac{\nu}{a_1^2\sqrt{\pi G \rho}}=0.01$. Viscosity destroys the gyroscopic stabilization of the Maclaurin spheroids on the interval $0.8127..<e<0.9529..$, which is stable
in the inviscid case \cite{C1969,C1984,RS1963}.}
\label{figcond2}
\end{figure}

The characteristic polynomial written for ${\mathbf L}_{v}(\lambda)$ yields the equation governing the growth rates of the ellipsoidal perturbations
in the presence of viscosity:
\be{grate}
25\Omega^2\mu^2+({\mathrm Re}\lambda+5\mu)^2 (\Omega^2-{\mathrm Re}\lambda^2-10{\mathrm Re} \lambda \mu-4b)=0.
\ee
The right panel of Figure~\ref{figcond2} shows that the growth rates \rf{grate} become positive beyond the Meyer-Liouville point.
Indeed, assuming ${\mathrm Re}\lambda=0$ in \rf{grate}, we reduce it to
$
50\mu^2(\Omega^2-2b)=0,
$
meaning that the growth rate vanishes when $\Omega^2=2b$ no matter how small the viscosity coefficient $\mu$ is. But, as we already know,
the equation $\Omega^2(e)=2b(e)$ determines exactly the Meyer-Liouville point, $e\approx 0.8127$.

      \begin{figure}
    \begin{center}
    \includegraphics[angle=0, width=0.32\textwidth]{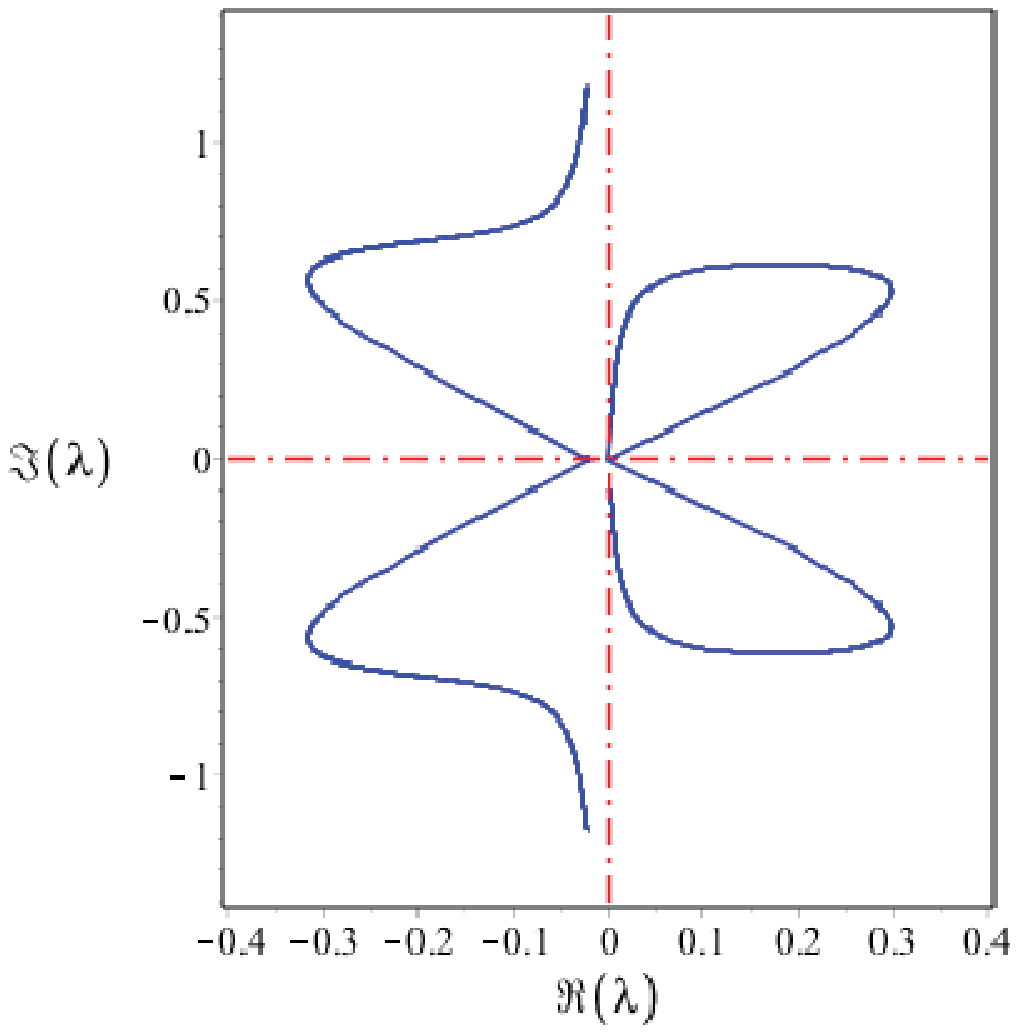}
\includegraphics[angle=0, width=0.32\textwidth]{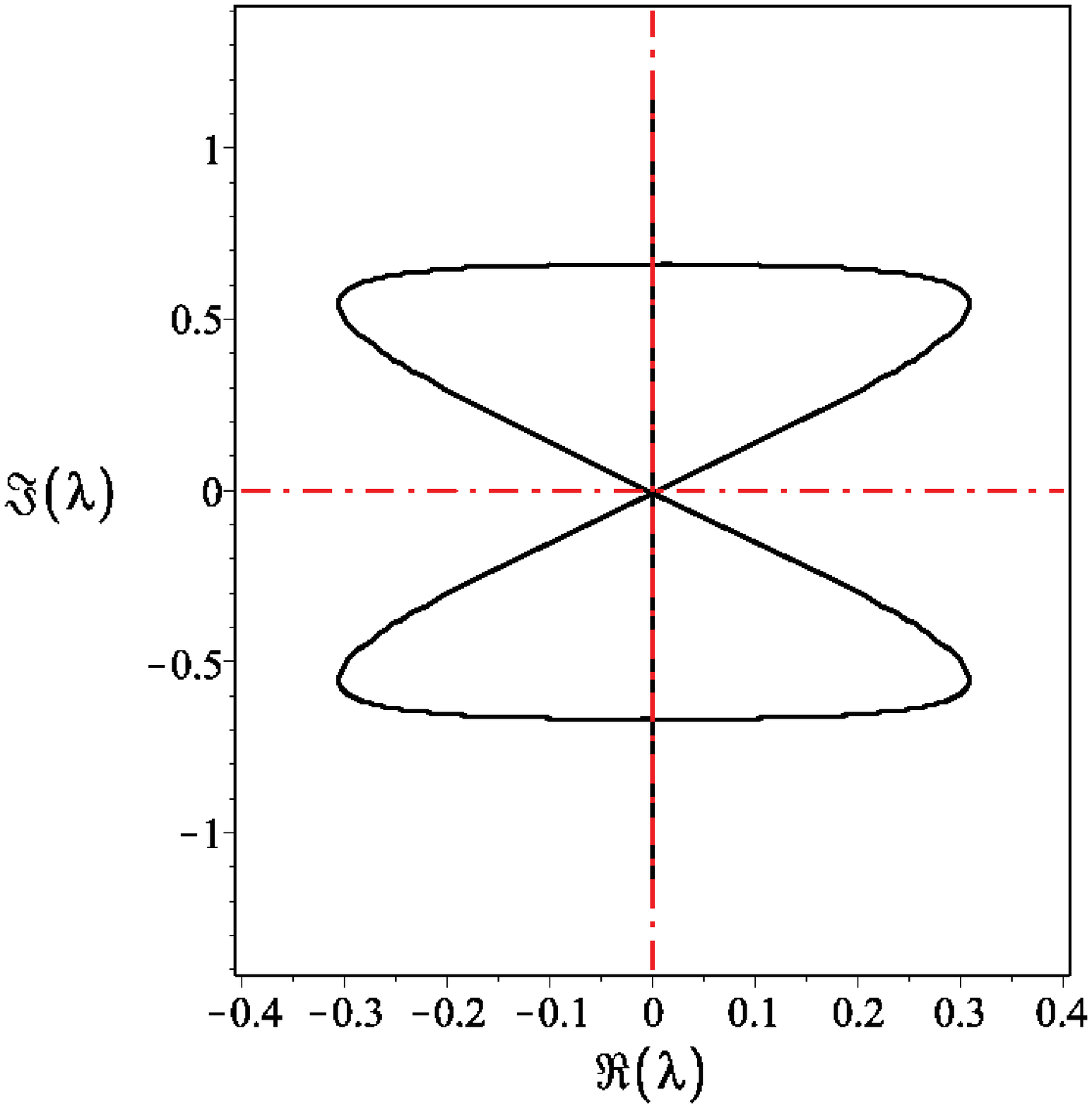}
\includegraphics[angle=0, width=0.32\textwidth]{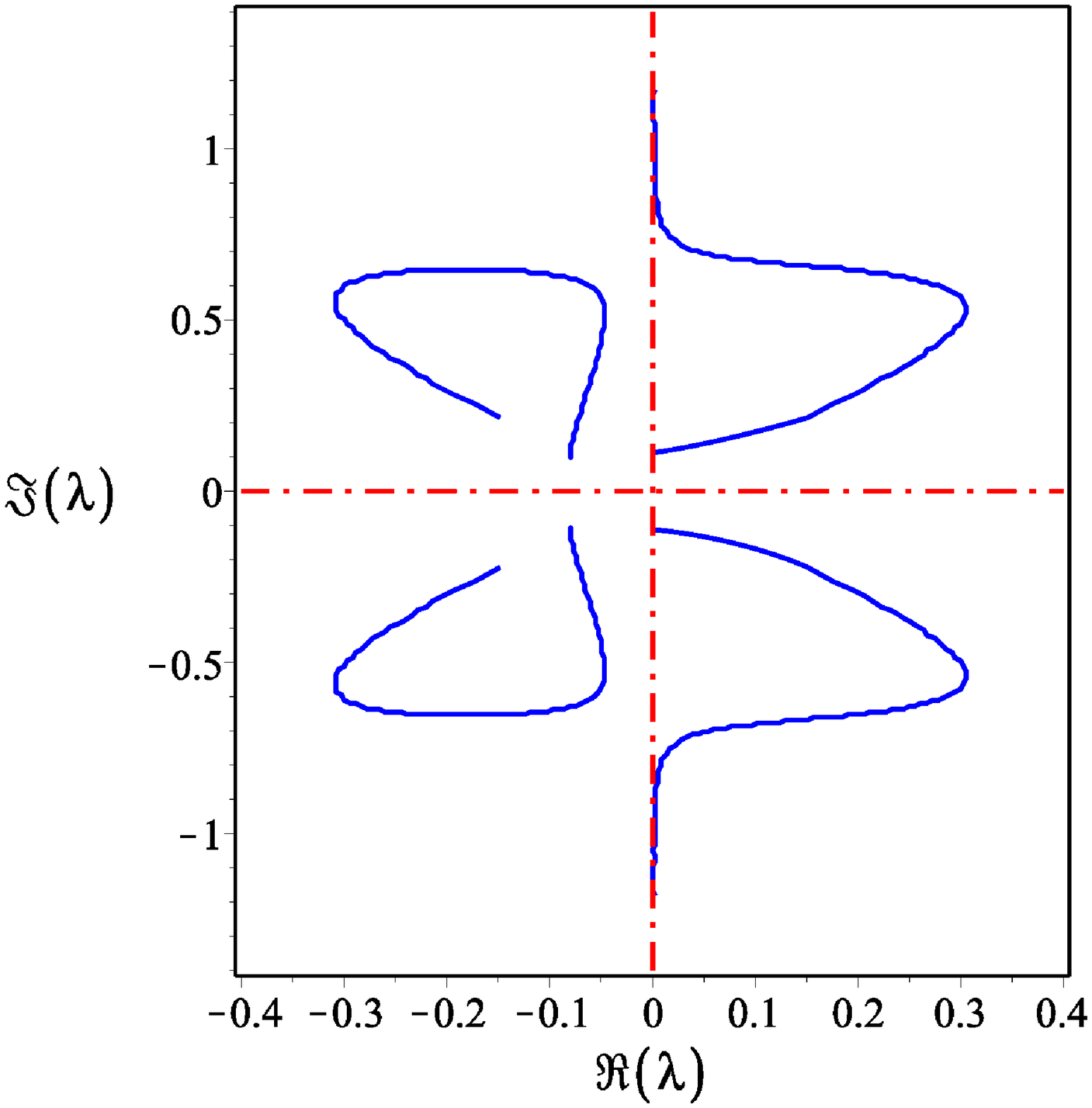}
    \end{center}
    \caption{Paths of the eigenvalues in the complex plane for (left) viscous Maclaurin spheroids with $\mu=\frac{\nu}{a_1^2\sqrt{\pi G \rho}}=0.002$,
(centre) Maclaurin spheroids without dissipation, and (right) inviscid Maclaurin spheroids with radiative losses for $\delta=0.05$. The
collision of two modes of the non-dissipative Hamiltonian system shown in the centre occurs at the Riemann critical value $e\approx 0.9529$.
Both types of dissipation destroy the collision and destabilize one of the two interacting modes at the Meyer-Liouville critical value $e \approx 0.8127$.}
    \label{figcond3}
    \end{figure}

It turns out, that the critical eccentricity of the viscous Maclaurin spheroid is equal to the
Meyer-Liouville value, $e\approx 0.8127$, even in the limit of vanishing viscosity, $\mu \rightarrow 0$, and thus does not converge to
the inviscid Riemann value $e\approx 0.9529$.

Viscous dissipation destroys the interaction of two modes at the Riemann critical point and destabilizes one of them beyond
the Meyer-Liouville point, showing an avoided crossing in the complex plane,  Figure~\ref{figcond3}(left).

Kelvin and Tait \cite{TT} hypothesised that the instability, which is stimulated by the presence of viscosity in the fluid,
will result in a slow, or
\textit{secular}, departure of the system from the unperturbed equilibrium of the Maclaurin family at the Meyer-Liouville point
and subsequent evolution along the Jacobi family, as long as the latter is
stable \cite{C1969}.

Therefore, a rotating, self-gravitating fluid mass, initially symmetric about the axis of rotation, can undergo
an axisymmetric evolution in which it first loses stability to a nonaxisymmetric disturbance,
and continues evolving along a non-axisymmetric family
toward greater departure from axial symmetry; then it undergoes a further loss of stability to a
disturbance tending toward splitting into two parts. Rigorous mathematical treatment of the validity of the fission theory of binary stars by Lyapunov and Poincar\'e has laid a foundation to modern nonlinear analysis.
In particular, it has led Lyapunov to the formulation of a general theory of stability of motion \cite{BKM2009}.

      \begin{figure}
    \begin{center}
    \includegraphics[angle=0, width=0.9\textwidth]{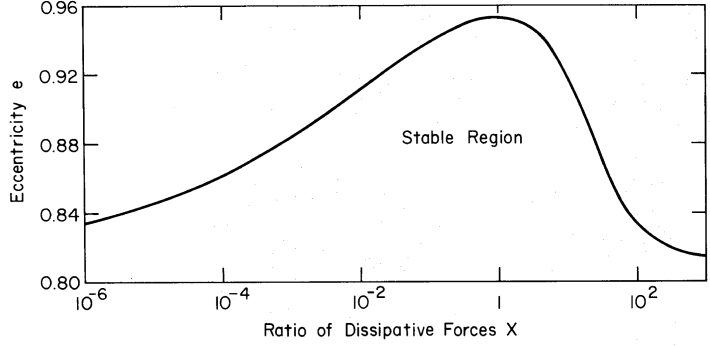}
    \end{center}
    \caption{Critical eccentricity in the limit of vanishing dissipation depends
    on the damping ratio, $X$, and attains its maximum (Riemann) value,
     $e\approx 0.9529$ exactly at $X=1$. As $X$ tends to zero or infinity, the critical value tends to the Meyer-Liouville value
$e\approx 0.8127$, \cite{LD1977,C1984}.}
    \label{figcond4}
    \end{figure}

In 1970 Chandrasekhar \cite{C1970} demonstrated that there exists another mechanism making Maclaurin spheroids
unstable beyond the Meyer-Liouville point of bifurcation, namely, the radiative losses due to emission of gravitational waves.
However, the mode that is made unstable by the radiation reaction is not the same one that is made unstable
by the viscosity, Figure~\ref{figcond3}(right).

In the case of the radiative damping mechanism stability is determined by the spectrum of the following matrix polynomial \cite{C1970}
$$
\textbf{L}_g(\lambda)=\lambda^2+\lambda(\textbf{G}+\textbf{D})+\textbf{K}+\textbf{N}
$$
that contains the matrices of gyroscopic, $\textbf{G}$, damping, $\textbf{D}$, potential, $\textbf{K}$, and
nonconservative positional, $\textbf{N}$, forces
$$
\textbf{G}=\frac{5}{2}\left(
                        \begin{array}{rr}
                          0 & -\Omega \\
                          \Omega & 0 \\
                        \end{array}
                      \right),\quad \textbf{D}=\left(
                                                 \begin{array}{cc}
                                                   \delta 16 \Omega^2(6b-\Omega^2) & -3\Omega/2 \\
                                                   -3\Omega/2 & \delta 16 \Omega^2(6b-\Omega^2) \\
                                                 \end{array}
                                               \right)
$$
$$
\textbf{K}=\left(
             \begin{array}{cc}
               4b-\Omega^2 & 0 \\
               0 & 4b-\Omega^2 \\
             \end{array}
           \right), \quad \textbf{N}=\delta\left(
                                       \begin{array}{cc}
                                         2q_1 & 2q_2 \\
                                         -q_2/2 & 2q_1 \\
                                       \end{array}
                                     \right),
                                     $$
where $\Omega(e)$ and $b(e)$ are given by equations \rf{maclaurin} and \rf{macb}.
Explicit expressions for $q_1$ and $q_2$ can be found in \cite{C1970}. The coefficient
$\delta=\frac{GMa_1^2(\pi G\rho)^{3/2}}{5c^5}$
is related to gravitational radiation reaction,
$G$ is the universal gravitation constant, $\rho$ the density of the fluid, $M$ the mass of the ellipsoid, and $c$ the velocity of light in vacuum.

In 1977 Lindblom and Detweiler \cite{LD1977} studied the combined effects of the gravitational radiation reaction and viscosity on the stability of the
Maclaurin spheroids. As we know, each of these dissipative effects induces a secular
instability in the Maclaurin sequence past the Meyer-Liouville point of bifurcation.
However, when both effects are considered together, the sequence of stable Maclaurin spheroids reaches past
the bifurcation point to a new point determined by the \textit{ratio}, $X=\frac{25}{2\Omega_0^4}\frac{\mu}{\delta}$, of the strengths of the viscous and radiation reaction forces, where $\Omega_0=0.663490...$.

Figure~\ref{figcond4} shows the critical eccentricity as a function of the damping ratio in the limit
of vanishing dissipation. This limit coincides with the inviscid Riemann point only at a particular damping ratio. At any other ratio,
the critical value is below the Riemann one and tends to the Meyer-Liouville value as this ratio tends either to zero or infinity.
Lindblom and Detweiler \cite{LD1977} correctly attributed the cancellation of the secular instabilities to the fact that viscous dissipation and radiation reaction cause
different modes to become unstable, see Figure~\ref{figcond3}. In fact, the mode destabilized by the fluid viscosity is the prograde moving spherical harmonic that appears to be retrograde in the frame
rotating with the fluid mass and the mode destabilized by the radiative losses is the retrograde moving spherical harmonic when it
appears to be prograde in the inertial frame \cite{A2003}. It is known \cite{Lamb1908,N1976} that to excite the positive energy mode
one must provide additional energy to the mode, while to excite the negative energy mode one must extract energy from the mode \cite{C95,Sa08}. The latter can be done by dissipation and the former by the non-conservative positional forces \cite{K07b, K09, K2013dg}. Both are present in the model by Lindblom and Detweiler \cite{LD1977}.

\begin{figure}
\begin{center}
\includegraphics[angle=0, width=0.87\textwidth]{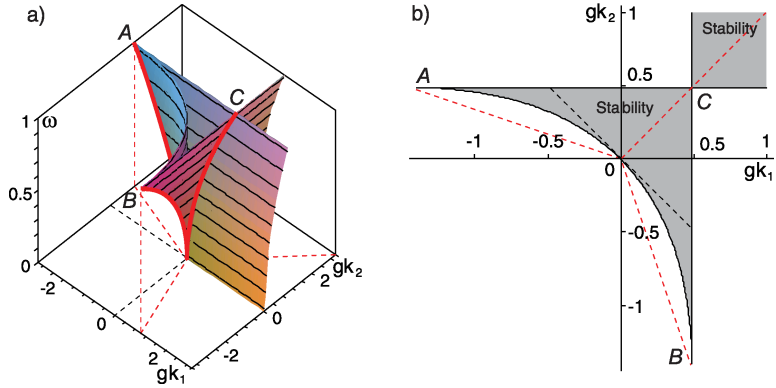}
\end{center}
\caption{(a) Stability domain \rf{bs} of the Brouwer's rotating vessel  and (b) its cross-section at $\omega=0.7$ \cite{K2011b,K2013dg}.}
\label{sadbro}
\end{figure}

\subsection{Brouwer's rotating vessel}

In 1918, Brouwer \cite{B18} published a simple model for the motion of a point mass near the bottom of a rotating vessel.
For an English translation see \cite{B76} pp. 665-675.
The shape of the vessel is described by a surface $S$, rotating  around a vertical axis with
constant angular velocity $\omega$. With the equilibrium chosen on the vertical axis at $(x, y)=(0, 0)$
on $S$, the linearized equations of motion without dissipation are
\begin{eqnarray}
\label{eqs1}
\ddot{x} - 2 \omega \dot{y} + (gk_1 - \omega^2) x &=& 0,\nn \\
\ddot{y} + 2 \omega \dot{x} + (gk_2 - \omega^2) y &=& 0.
\end{eqnarray}
The constants $k_1$ and $k_2$ are the $x, y$-curvatures of $S$ at $(0, 0)$, $g$ is the gravitational constant.
Suppose that $k_1 \geq k_2$.\\

\noindent
{\bf Stability and instability without friction}\\
Assuming there is no damping, there are the following three cases:
\begin{itemize}
\item $0<k_2<k_1$ (single-well at equilibrium).\\
Stability iff $0< \omega^2< gk_2$ (slow rotation)
or $\omega^2> gk_1$ (fast rotation).
\item $k_2<0$ and $k_1>0, \, k_1> -k_2$ (saddle at equilibrium).\\
Stability iff  $\omega^2> gk_1$.
\item $k_2<0$ and $k_1>0, \, k_1< -k_2$ (saddle).\\
If $3k_1+k_2 <0$: instability.\\
If $3k_1+k_2 >0$: stability if
\be{bs} gk_1 < \omega^2 < - \frac{g}{8} \frac{(k_1-k_2)^2}{(k_1+k_2)}.\ee
\item $k_1<0, \, k_2<0$: instability.
\end{itemize}

\noindent
{\bf Stability of triangular libration points $L_4$ and $L_5$}\\
Brouwer's rotating vessel model includes both a well with two positive curvatures $k_1$ and $k_2$ and a saddle with the curvatures of opposite signs. Remarkably, the latter case describes stability of triangular libration points $L_4$ and $L_5$ (discovered by Lagrange in 1772) in the restricted circular three-body problem of celestial mechanics. Indeed, the linearized equation for this problem is \rf{eqs1} where $\omega=1$,
\ba{L4b}
gk_1&=&-\frac{1}{2}+\frac{3}{2}\sqrt{1-3\mu(1-\mu)}\ge \frac{1}{4},\nn\\
gk_2&=&-\frac{1}{2}-\frac{3}{2}\sqrt{1-3\mu(1-\mu)}\le -\frac{5}{4},
\ea
$\mu=\frac{m_1}{m_1+m_2}$, and $m_1$ and $m_2$ are the masses of the two most massive bodies (in comparison with each of them the mass of the third body is assumed to be negligible) \cite{A1970,G1843,R2002}.
Since $k_1$ and $k_2$ are of opposite signs for all $0\le\mu\le 1$, the linear stability of the triangular Lagrange equilibriums is determined by the stability conditions \rf{bs} for the rotating saddle \cite{BB1994}. Note that the coefficients \rf{L4b} satisfy the constraint $gk_1+gk_2=-1$.  Intervals of intersection of this line with the narrow corners of the curvilinear triangle in Figure~\ref{sadbro}b correspond to the stable Lagrange points. After substitution of \rf{L4b} into \rf{sadbro}, we reproduce the classical condition for their linear stability first established by Gascheau in 1843 \cite{G1843}
$$
(gk_1-gk_2)^2-8=1-27\mu(1-\mu)>0.
$$
Stability in a rotating saddle potential is a subject of current active discussion in respect with the particle trapping \cite{BB1994,BC2018,K2011b,JMT2017}. An effect of slow precession of trajectories of the trapped particles in a rotating saddle potential \cite{KL2017,S1998} has recently inspired new works leading to the improvement of traditional averaging methods \cite{CL2018,YH2018}.\\

\begin{figure}
\begin{center}
\includegraphics[angle=0, width=0.4\textwidth]{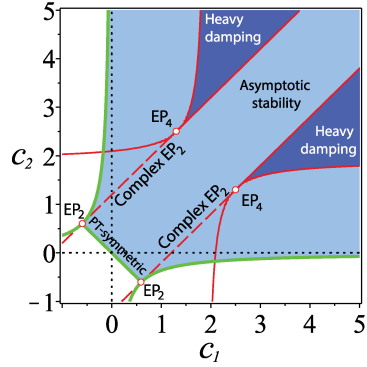}
\end{center}
\caption{Domain of asymptotic stability of the damped Brouwer's rotating vessel in the plane of the damping coefficients for $gk_1=1$, $gk_2=1$, and $\omega=0.3$ \cite{K2013dg,K2013pt}.}
\label{sadbrod}
\end{figure}

\noindent
{\bf Destabilization by friction}\\
Adding constant (Coulomb) damping, Brouwer \cite{B18} finds a number of instability cases.
The equations of motion become in this case:
\begin{eqnarray} \label{eqs1+d}
\ddot{x} - 2 \omega \dot{y} +c_1 \dot{x}+ (gk_1 - \omega^2)x &=& 0,\\
\ddot{y} + 2 \omega \dot{x} +c_2 \dot{y} + (gk_2 - \omega^2)y &=& 0.
\end{eqnarray}
The friction constants $c_1, c_2$ are positive. The
characteristic (eigenvalue) equation takes the form:
\[ \lambda^4 + a_1 \lambda^3+ a_2 \lambda^2 +a_3 \lambda +a_4=0, \]
with $a_1= c_1+c_2$, \\
$a_2= g(k_1+k_2)+ 2 \omega^2 + c_1c_2$\\
$a_3= c_1(gk_2- \omega^2) + c_2(gk_1- \omega^2)$, \\
$a_4=(gk_1- \omega^2)(gk_2- \omega^2)$.

There are two cases that drastically change the stability (see also Bottema \cite{Bo76}) :
\begin{itemize}
\item $0<k_2 < k_1$ (single-well).\\
Stability iff $0< \omega^2< gk_2$.\\
The fast rotation branch $\omega^2> gk_1$ has vanished.

\item $k_2<0$ and $k_1>0, \, k_1< -k_2$ (saddle).\\
The requirement $a_4>0$ produces $0<gk_1< \omega²$.
This is not compatible with $a_3<0$ so  a saddle is always unstable
with any size of positive damping.
\end{itemize}
Brouwer studied this model probably to use in his lectures. In a correspondence with O. Blumenthal and
G. Hamel he asked whether the results of the calculations were known; see \cite{B76} pp. 677-686. Hamel confirmed that the results were correct and surprising, but there is no reference to older literature in this correspondence. See also \cite{Bo76}.\\

\noindent
{\bf Indefinite damping and PT-symmetry}\\
Brouwer's model for the case of a rotating well with two positive curvatures has a direct relation to rotordynamics as it contains the Jeffcott rotor model, see for instance \cite{K2011b}. Dissipation induces instabilities in this model at high speeds of rotation $\omega$ \cite{Bo63,C95,S33}, of course, under the assumption that the damping coefficients are both positive. But what happens if we relax this constraint and extend the space of damping parameters to negative values? It turns out that at low speeds $\omega$ the domain of asymptotic stability spreads to the area of negative damping, Figure~\ref{sadbrod}. Even more,
a part of the neutral stability curve belongs to the line $c_1+c_2=0$ where one of the damping coefficients is positive and the other one is negative. On this line the system is invariant under time and parity reversion transformation and is therefore $PT$-symmetric. Its eigenvalues are imaginary in spite of the presence of the loss (positive damping) and gain (negative damping) in the system. $PT$-symmetric systems with the indefinite damping can easily be realized in the laboratory experiments \cite{Kottos2011}. Recent study \cite{Z2018} shows their connection to complex G-Hamiltonian systems \cite{KI2017,YS75}.

 \section{Ziegler's paradox}
  We already know that Greenhill's analysis of stability of the Kelvin gyrostat \cite{G1880} inspired the works of Sobolev \cite{S2006} and Pontryagin \cite{P44} which have led to the development of the theory of spaces with indefinite inner product \cite{Bognar1974,K1983}. Another work of Greenhill \cite{G1883} ultimately brought about the famous Ziegler paradox \cite{Z52}. As Gladwell remarked in his historical account of the genesis of the field of nonconservative stability \cite{G1990}, ``It was Greenhill who started the trouble though he never knew it.''

  Motivated by the problem of buckling of propeller-shafts of steamers Greenhill (1883) analyzed stability of an elastic shaft of a circular cross-section under the action of a compressive force and an axial torque \cite{G1883}. He managed to find the critical torque that causes buckling of the shaft for a number of boundary conditions. For the clamped-free and the clamped-hinged shaft loaded by an axial torque the question remained open until Nicolai in 1928 reconsidered a variant of the clamped-hinged boundary conditions, in which the axial torque is replaced with the \textit{follower} torque \cite{N1928}. The vector of the latter is directed along the tangent to the deformed axis of the shaft at the end point \cite{G1990}.

  Nicolai had established that no nontrivial equilibrium configuration of the shaft exists different from the rectilinear one,
meaning stability for all magnitudes of the follower torque. Being unsatisfied with this
overoptimistic result, he assumed that the equilibrium method does not work properly in the case of the follower torque.
He decided to study small oscillations of the shaft about its rectilinear configuration using what is now known as the
Lyapunov stability theory \cite{L1992} that, in particular, can predict instability via eigenvalues of the linearized problem.

Surprisingly, it turned out that there exist eigenvalues with positive real parts (instability) for all magnitudes of the torque,
meaning that the critical value of the follower torque for an elastic shaft of a circular cross-section is actually zero.
Because of its unusual behavior, this instability phenomenon received a name ``Nicolai's paradox'' \cite{G1990,LFA2016,N1928}.

In 1951-56 Hans Ziegler of the ETH Z\"urich re-considered the five original Greenhill problems with the Lyapunov approach and found that the shaft is unstable in cases of the clamped-free and the clamped-hinged boundary conditions for all values of the axial torque, just as in Nicolai's problem with the follower torque \cite{Z53}. Moreover, the non-self-adjoint boundary eigenvalue problem for the Greenhill's shaft with the axial torque turned out to be a Hermitian adjoint of the non-self-adjoint boundary eigenvalue problem for the Greenhill's shaft with the follower torque \cite{Bo63}.

 \begin{figure}
 \begin{center}
\includegraphics[width=0.8\textwidth]{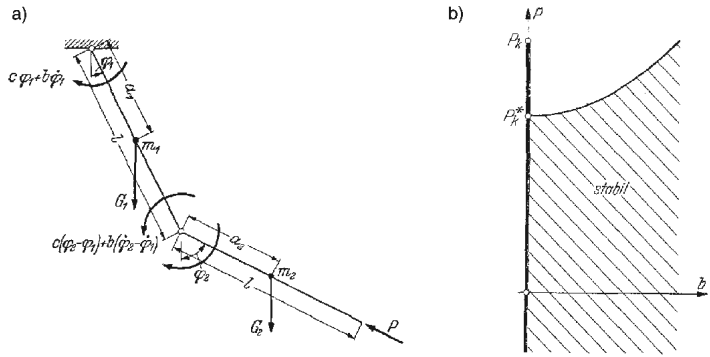}
  \end{center}
\caption{\label{fig1} Original drawings from Ziegler's work of 1952 \cite{Z52}: (a) double linked pendulum under the follower load $P$, (b) (bold line) stability interval of the undamped pendulum and (shaded area) the domain of asymptotic stability of the damped pendulum with equal coefficients of dissipation in both joints.
 If $b=0$ we have no dissipation and the system is marginally stable for $P < P_k$.}
\label{dpend}
\end{figure}

 In 1952, inspired by the paradoxes of Greenhill-Nicolai follower torque problems, Ziegler introduced the notion of the {\em follower force} and published a paper \cite{Z52} that became widely known in the engineering community, in particular among those interested in theoretical mechanics. It was followed by a second paper \cite{Z53} that added more details.

 Ziegler considered a double pendulum consisting of two rigid rods of length $l$ each.  The
 pendulum is attached at one of the endpoints and can swing freely in a vertical plane; see Figure~\ref{dpend}. The
 angular deflections with respect to the vertical are denoted by $\phi_1, \phi_2$, two masses $m_1$ and $m_2$
 resulting in the external forces $G_1$ and $G_2$ are concentrated at the distances $a_1$ and $a_2$ from the joints. At the joints
 we have elastic restoring forces of the form $c \phi_1$, $c(\phi_2- \phi_1)$ and internal damping torques
 \[  b_1 \frac{d \phi_1}{dt},\quad b_2\left(\frac{d \phi_2}{dt} - \frac{d \phi_1}{dt}\right). \]
 So if $b_1=b_2=0$ we have no dissipation.
 We impose a follower force $P$ on the lowest hanging rod, see Figure~\ref{dpend}. We consider only the quadratic
 terms of kinetic and potential energy. With these assumptions
  the kinetic energy $T$ of the system is:
 \begin{equation}
 T= \frac{1}{2} \left[(m_1a_1^2+m_2l^2)\dot{\phi_1}^2 + 2m_2la_2 \dot{\phi_1} \dot{\phi_2}+m_2 a_2^2 \dot{\phi_2}^2 \right].
 \end{equation}
 A dot denotes differentiation with respect to time $t$. The potential energy $V$ reads:
 \begin{equation}
 V = \frac{1}{2}\left[(G_1a_1+G_2l+2c)\phi_1^2- 2c \phi_1 \phi_2+ (G_2a_2+c)\phi_2^2 \right].
 \end{equation}
 This leads to the generalised dissipative and non-conservative forces $Q_1, Q_2$:
 \begin{equation}
 Q_1= Pl(\phi_1- \phi_2)-((b_1+b_2)\dot{\phi}_1-b_2 \dot{\phi}_2),\,\,Q_2= b_2(\dot{\phi}_1-  \dot{\phi}_2).
 \end{equation}

Writing the Lagrange's equations of motion $\dot L_{\dot\varphi_i}-L_{\varphi_i}=Q_i$, where $L=T-V$ and
 assuming $G_1=0$ and $G_2=0$ for simplicity, we find
\ba{z4}
\left(
  \begin{array}{rr}
    m_1 a_1^2+m_2l^2 & m_2 l a_2 \\
    m_2la_2 & m_2 a_2^2 \\
  \end{array}
\right)\left(
         \begin{array}{c}
           \ddot \varphi_1 \\
           \ddot \varphi_2 \\
         \end{array}
       \right)&+&
       \left(
  \begin{array}{rr}
    b_1+b_2 & -b_2 \\
    -b_2 & b_2 \\
  \end{array}
\right)\left(
         \begin{array}{c}
           \dot \varphi_1 \\
           \dot \varphi_2 \\
         \end{array}
       \right)\nn\\&+&
              \left(
  \begin{array}{rr}
    -Pl+2c & Pl-c \\
    -c & c \\
  \end{array}
\right)\left(
         \begin{array}{c}
           \varphi_1 \\
           \varphi_2 \\
         \end{array}
       \right)=0.
\ea
The stability analysis of equilibrium follows the standard procedure.
With the substitution  $\varphi_i=A_i \exp(\lambda t)$, equation \rf{z4} yields a
$4$-dimensional eigenvalue problem with respect to the spectral parameter $\lambda$.

Putting $m_1=2m$, $m_2=m$, $a_1=a_2=l$, $b_1=b_2=b$ and assuming that internal damping is absent $(b=0)$, Ziegler found from the characteristic equation that the vertical equilibrium position of the pendulum looses its stability when the magnitude of the follower force exceeds the critical value $P_k$, where
\be{z5}
P_k=\left(\frac{7}{2}-\sqrt{2} \right)\frac{c}{l}\simeq2.086\frac{c}{l}.
\ee

\begin{figure}
\begin{center}
\includegraphics[width=0.9\textwidth]{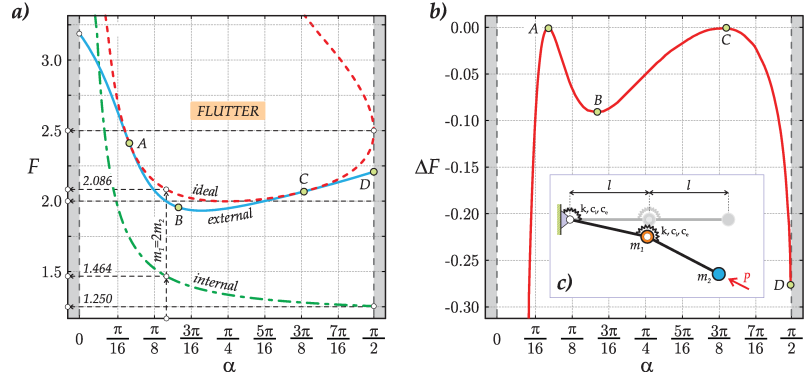}
\end{center}
\caption{(a) The (dimensionless) follower force $F$, shown as a function of the (transformed via $\cot{\alpha}=m_1/m_2$) mass ratio $\alpha$, represents the flutter domain of (dashed/red line) the undamped, or \lq ideal', Ziegler pendulum \cite{SW1975} and the flutter boundary of the dissipative system in the limit of vanishing (dot-dashed/green line) internal and (continuous/blue line) external damping. (b) Discrepancy $\Delta F$ between the critical flutter load for the ideal Ziegler pendulum and for the same structure calculated in the limit of vanishing external damping. The discrepancy quantifies the Ziegler's paradox due to external (air drag) damping \cite{T2016}.}
\label{pflu}
\end{figure}

In the presence of damping $(b > 0)$ the Routh-Hurwitz condition yields the new critical follower load
that depends on the square of the damping coefficient $b$
\be{z6}
P_k(b)=\frac{41}{28}\frac{c}{l}+\frac{1}{2}\frac{b^2}{ml^3}.
\ee
Ziegler found that the domain of asymptotic stability for the damped pendulum is given by the inequalities
$P<P_k(b)$ and $b>0$ and he plotted it against the stability interval $P<P_k$ of the undamped system, Figure~\ref{fig1}(b). Surprisingly, the limit of the critical load $P_k(b)$ when $b$ tends to zero turned out to be significantly lower than the critical load of the undamped system
\be{z7}
P_k^*=\lim_{b\rightarrow0}P_k(b)=\frac{41}{28}\frac{c}{l}\simeq1.464\frac{c}{l}<P_k.
\ee

Note that in the original work of Ziegler, formula \rf{z6} contains a misprint which yields linear dependency of the critical follower load on the damping coefficient $b$. Nevertheless, the domain of asymptotic stability found in \cite{Z52} and reproduced in Figure~\ref{fig1}(b), is correct.

Ziegler limited his original calculation to a particular mass distribution, $m_1=2m_2$, and took into account only internal damping in the joints, neglecting, e.g, the air drag (an external damping). Later studies confirmed that the Ziegler paradox is a generic phenomenon and exists at all mass distributions, both for internal and external damping \cite{BK2018,B2018,SLK2019,T2016}, see Figure~\ref{pflu}.

After invention of robust methods of practical realization of follower forces \cite{BN2011} the Ziegler paradox was immediately observed in the recent laboratory experiments \cite{BK2018,B2018}. Nowadays follower forces find new applications in cytosceletal dynamics \cite{BD2016,DC2017} and acoustics of friction \cite{HG03}. In general, the interest to mathematical models involving nonconservative positional forces (known also as circulatory \cite{Z53} or curl \cite{BS2016} forces)
is growing both in traditional areas such as energy harvesting and fluid-structure interactions \cite{HD92,PMB2017} and in rapidly emerging new research fields of optomechanics \cite{SD2017} and light robotics \cite{LR2017}.

 \section{Bottema's analysis of Ziegler's paradox}
In 1956, in the journal `Indagationes Mathematicae', there appeared an article by Oene Bottema (1901-1992) \cite{B56},
then Rector Magnificus of the Technical University of Delft and an expert in classical geometry and mechanics,
that outstripped later findings for decades.
Bottema's work on stability in 1955 \cite{B55} can be seen as an introduction, it was directly motivated by
Ziegler's paradox. However, instead of examining the particular model of Ziegler, he studied in \cite{B56} a much
more general class of non-conservative systems.

Following \cite{B55,B56}, we consider a holonomic scleronomic linear system with two degrees of
freedom, of which the coordinates $x$ and $y$ are chosen in such a way that the kinetic energy is
$
T=({\dot x}^2+{\dot y}^2)/2.
$
Hence the Lagrange equations of small oscillations
near the equilibrium configuration $x = y = 0$ are as follows
\ba{b2}
\ddot x +a_{11}  x +a_{12}  y +b_{11} \dot x +b_{12} \dot y&=& 0,\nn \\
\ddot y +a_{21}  x +a_{22}  y +b_{21} \dot x +b_{22} \dot y&=& 0,
\ea
where $a_{ij}$ and $b_{ij}$ are constants, ${\bf A}:=(a_{ij})$ is the matrix of the forces depending on
the coordinates, ${\bf B}:=(b_{ij})$ of those depending on the velocities. If $\bf A$ is
symmetrical and while disregarding the damping associated with
the matrix $\bf B$, there exists a potential energy function
$V = (a_{11} x^2 + 2 a_{12}x y + a_{22} y^2)/2$,
if it is antisymmetrical, the forces are circulatory. When the matrix $\bf B$ is
symmetrical, we have a non-gyroscopic damping force, which is positive when
the dissipative function
$(b_{11} x^2 + 2 b_{12}x y + b_{22} y^2)/2$
is positive definite. If $\bf B$ is antisymmetrical the forces depending on the
velocities are purely gyroscopic.

The matrices $\bf A$ and $\bf B$ can both be written uniquely as the sum of symmetrical and
antisymmetrical parts: ${\bf A}={\bf K}+{\bf N}$ and ${\bf B}={\bf D}+{\bf G}$, where
\be{b3}
{\bf K}=\left(
          \begin{array}{cc}
            k_{11} & k_{12} \\
            k_{21} & k_{22} \\
          \end{array}
        \right), \quad {\bf N}=\left(
                  \begin{array}{rr}
                    0 & \nu \\
                    -\nu & 0 \\
                  \end{array}
                \right),\quad
                {\bf D}=\left(
          \begin{array}{cc}
            d_{11} & d_{12} \\
            d_{21} & d_{22} \\
          \end{array}
        \right), \quad {\bf G}=\left(
                  \begin{array}{rr}
                    0 & \Omega \\
                    -\Omega & 0 \\
                  \end{array}
                \right),
\ee
with $k_{11}=a_{11}$, $k_{22}=a_{22}$, $k_{12}=k_{21}=(a_{12}+a_{21})/2$, $\nu=(a_{12}-a_{21})/2$ and $d_{11}=b_{11}$, $d_{22}=b_{22}$, $d_{12}=d_{21}=(b_{12}+b_{21})/2$, $\Omega=(b_{12}-b_{21})/2$.

Disregarding damping, the system \rf{b2} has a potential energy function  when $\nu= 0$, it is purely circulatory
for $k_{11} = k_{12}=k_{22} = 0$, it is non-gyroscopic for $\Omega= 0$, and has no damping when
$d_{11} = d_{12} = d_{22} = 0$. If damping exists, we suppose in this section that it is positive.

In order to solve the linear stability problem for equations \rf{b2} we put $x = C_1 \exp(\lambda t)$,
$y = C_2 \exp(\lambda t)$ and
obtain the characteristic equation for the frequencies of the small oscillations around equilibrium
\be{b5a}
Q:=\lambda^4 + a_1\lambda^3 + a_2\lambda^2 + a_3\lambda + a_4=0,
\ee
where \cite{K03a,K03b,K07a}
\be{b5}
a_1 = {\rm tr}{\bf D},\quad
a_2 = {\rm tr}{\bf K} + \det {\bf D} + \Omega^2,\quad
a_3 = {\rm tr}{\bf K}{\rm tr} {\bf D} - {\rm tr} {\bf KD} + 2\Omega\nu ,\quad
a_4 = \det{\bf K} + \nu^2.
\ee
For the equilibrium to be stable all roots of the characteristic equation \rf{b5a} must be semi-simple and have real parts which
are non-positive.

It is always possible to write, in at least one way, the left hand-side as
the product of two quadratic forms with real coefficients,
$
Q=(\lambda^2+p_1\lambda+q_1)(\lambda^2+p_2\lambda+q_2).
$
Hence
\be{b7a}
a_1=p_1+p_2,\quad a_2=p_1p_2+q_1+q_2,\quad a_3=p_1q_2+p_2q_1,\quad a_4=q_1q_2.
\ee

For all the roots of the equation \rf{b5a} to be in the left side of the complex plane $(L)$ it is obviously necessary and sufficient
that $p_i$ and $q_i$ are positive. Therefore in view of \rf{b7a} we have: a necessary
condition for the roots of $Q=0$ having negative real parts is $a_i>0$ $(i = 1, 2, 3,4)$.
This system of conditions however is not sufficient, as the example
$(\lambda^2-\lambda+2)(\lambda^2+2\lambda+3)=\lambda^4+\lambda^3+3\lambda^2+\lambda+6$ shows. But if $a_i>0$ it is not
possible that either one root of three roots lies in $L$ (for then $a_4\le0$); it
is also impossible that no root is in it (for, then $a_4 \le 0$). Hence if $a_i > 0$ at
least two roots are in $L$; the other ones are either both in $L$, or both on
the imaginary axis, or both in $R$. In order to distinguish between these cases we
deduce the condition for two roots being on the imaginary axis. If $\mu i$ ($\mu \ne 0$ is real)
is a root, then $\mu^4-a_2\mu^2+a_4=0$ and $-a_1\mu^2+a_3=0$.
Hence $H:=a_1^2a_4+a_3^2-a_1a_2a_3=0$.
Now by means of \rf{b7a} we have
\be{bb6}
H=-p_1p_2(a_1a_3+(q_1-q_2)^2).
\ee

In view of $a_1> 0$, $a_3>0$ the second factor is positive; furthermore
$a_1= p_1+p_2 > 0$, hence $p_1$ and $p_2$ cannot both be negative. Therefore $H < 0$
implies $p_1> 0$, $p_2> 0$, for $H = 0$ we have either $p_1 = 0$ or $p_2=0$ (and not both,
because $a_3>0$), for $H>0$ $p_1$ and $p_2$ have different signs. We see from
the decomposition of the polynomial \rf{b5a} that all its roots are in $L$ if $p_1$ and $p_2$ are positive.

Hence: a set of necessary and sufficient conditions for all roots of \rf{b5a}
to be on the left hand-side of the complex plane is
\be{bb6a}
a_i>0~~(i=1,2,3,4), \quad H<0.
\ee

We now proceed to the cases where all roots have non-positive
real parts, so that they lie either in $L$ or on the imaginary axis.
If three roots are in $L$ and one on the imaginary axis, this root must
be $\lambda = 0$. Reasoning along the same lines as before we find that necessary
and sufficient conditions for this are
$a_i>0~~(i=1,2,3)$, $a_4=0$, and $H<0$.
If two roots are in $L$ and two (different) roots on the imaginary axis we have
$p_1>0$, $q_1>0$, $p_2=0$, $q_2>0$ and the conditions are
$a_i>0~~(i=1,2,3,4)$ and $H=0$.
If one root is in $L$ and three are on the imaginary axis, then $p_1>0$, $q_1=0$, $p_2=0$,
$q_2>0$ and the conditions are
$a_i>0~~(i=1,2,3)$, $a_4=0$, and $H=0$.

The obtained conditions are border cases of \rf{bb6a}. This does not occur
with the last type we have to consider: all roots are on the imaginary axis. We
now have $p_1=0$, $p_2=0$, $q_1>0$, $q_2>0$. Hence $a_2 >0$, $a_4>0$, $a_1=a_3=0$ and
therefore $H = 0$. This set of relations is necessary, but not sufficient, as
the example $Q= \lambda^4 + 6\lambda^2+ 25 = 0$ (which has two roots in $L$ and two in the righthand side of the complex plane $(R)$)
shows. The proof given above is not valid because as seen from \rf{bb6a}, $H=0$
does not imply now $p_1p_2= 0$, the second factor being zero for $a_1a_3 = 0$
and $q_1 =q_2$. The condition can of course easily be given; the equation \rf{b5a} is
$\lambda^4+a_2\lambda^2+a_4=0$
and therefore it reads $a_2>0$, $a_4>0$, $a_2^2>4a_4$.

Summing up we have: all roots of \rf{b5a} (assumed to be different) have
non-positive real parts if and only if one of the two following sets of
conditions is satisfied \cite{B56}
\ba{bb11}
A:&~& a_1>0,~a_2>0,~a_3>0,~a_4\ge0,~a_2\ge \frac{a_1^2a_4+a_3^2}{a_1a_3},\nn \\
B:&~& a_1=0,~a_2>0,~a_3=0,~a_4>0,~a_2>2\sqrt{a_4}.
\ea

Note that $a_1$ represents the damping coefficients $b_{11}$ and $b_{22}$ in the system.
One could expect $B$ to be a limit of $A$, so that for $a_1 \rightarrow 0$, $a_3 \rightarrow 0$
the set $A$ would continuously tend to $B$. {\em That is not the case}.

Remark first of all that the roots of \rf{b5a} never lie outside $R$ if $a_1 = 0$, $a_3 \ne 0$ (or $a_1 \ne  0$,
$a_3=0$). Furthermore, if $A$ is satisfied and we take $a_1=\varepsilon b_1$, $a_3=\varepsilon b_3$, where
$b_1$ and $b_3$ are fixed and $\varepsilon \rightarrow 0$, the last condition of $A$ reads $(\varepsilon \ne 0)$
$$
a_2 > \frac{b_1^2a_4+b_3^2}{b_1b_3}=g_1
$$
while for $\varepsilon = 0$ we have
$$
a_2> 2 \sqrt{a_4}=g_2.
$$
Obviously we have \cite{B56}
$$
g_1-g_2=\frac{(b_1 \sqrt{a_4}-b_3)^2}{b_1b_3}
$$
so that $(g_1>g_2)$ but for $b_3=b_1\sqrt{a_4}$. That means that in all cases where $b_3\ne b_1\sqrt{a_4}$ we have a discontinuity in our
stability condition.

\begin{figure}
\includegraphics[width=0.9\textwidth]{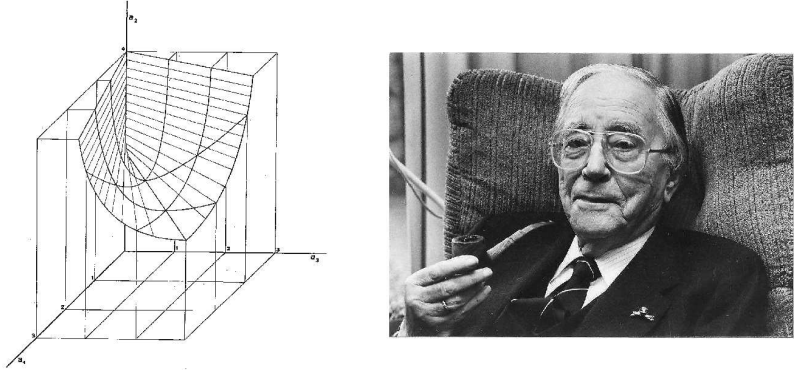}
\caption{\label{fig2} Original drawing (left) from the 1956 work \cite{B56} of Oene Bottema (right),
showing the domain of the asymptotic stability of the real polynomial of fourth order and of the two-dimensional non-conservative system
with Whitney's umbrella singularity discussed in the sequel. The ruled surface (called $V$ in the text) is given by equation
(\ref{bb7}).}
\end{figure}

In 1987, Leipholz remarked in his monograph on stability theory \cite{L1987} that ``Independent works of Bottema and Bolotin for \textit{second-order systems} have shown that in the non-conservative case and for different damping coefficients the stability condition is discontinuous with respect to the undamped case.'' However, Leipholz did not mention that, in contrast to Bolotin, Bottema illustrated the phenomenon of the discontinuity in a remarkable \textit{geometric diagram}, first published in \cite{B56} and reproduced in Figure~\ref{fig2}.

Following Bottema \cite{B56} we substitute in \rf{b5a} $\lambda=c\mu$, where $c$ is the positive fourth root
of $a_4>0$. The new equation reads $P:=\mu^4+ b_1\mu^3+ b_2\mu^2+b_3\mu+1 = 0$,
where $b_i=a_i/c^i$ $(i= 1,2,3,4)$. If we substitute $a_i=c^ib_i$ in $A$ and $B$ we get
the same condition as when we write $b_i$ for $a_i$, which was to be expected,
because if the roots of \rf{b5a} are outside $R$, those of $P=0$ are also outside $R$
and inversely. We can therefore restrict ourselves to the case $a_4 = 1$, so
that we have only three parameters $a_1$, $a_2$, $a_3$. We take them as coordinates
in an orthogonal coordinate system.

The condition $H=0$ or
\be{bb7}
a_1a_2a_3 = a_1^2 +a_3^2
\ee
is the equation of a surface $V$ of the third degree, which we have to
consider for $a_1 \ge 0$, $a_3 \ge 0$, Figure~\ref{fig2}. Obviously $V$ is a \textit{ruled surface}, the line $a_3=ma_1$,
$a_2=m+1/m$ $(0<m<\infty)$ being on $V$. The line is parallel to the $0a_1a_3$-plane
and intersects the $a_2$-axis in $a_1=a_3= 0$, $a_2=m+ 1/m\ge 2$. The $a_2$-axis
is the double line of $V$, $a_2>2$ being its active part.
Two generators pass through each point of it; they coincide for $a_2 = 2$ $(m = 1)$, and for $a_2\rightarrow\infty$
their directions tend to those of the $a_1$ and $a_3$-axis $(m=0, m=\infty)$.
The conditions $A$ and $B$ express that the image point $(a_1, a_2, a_3)$ lies
on $V$ or above $V$. The point $(0, 2, 0)$ is on $V$, but if we go to the $a_2$-axis
along the line $a_3=ma_1$ the coordinate $a_2$ has the limit $m+ 1/m$, which is
$>2$ but for $m= 1$.

Note that we started off with 8 parameters in Eq. \rf{b2}, but that the surface $V$ bounding the stability domain
is described by 3 parameters. It is described by a map of $E^2$ into $E^3$ as in Whitney's papers \cite{W43,W44}. Explicitly, a
transformation of (19) to (2) is given by
$$
a_1 =\frac{1}{2}y_3 + w,~~ a_2 = 2 + y_2,~~ a_3 = -\frac{1}{2}y_3 + w
$$
with $w^2 = \frac{1}{4}y_3^2+ y_1y_2$.

Returning to the non-conservative system \rf{b2} $(\nu \ne 0)$, with damping,
but without gyroscopic forces, so $\Omega = 0$, and assuming as in \cite{B55} that $k_{12}=0$, $k_{11}> 0$, and $k_{22} > 0$ (a similar setting but with $d_{12}=0$ and $k_{12}\ne0$ was considered later by Bolotin in \cite{Bo63}), we find that the condition for stability $H\le0$ reads
\ba{bb2}
\nu^2 &<& \frac{(k_{11}-k_{22})^2}{4}\\
&-&\frac{(d_{11}-d_{22})^2(k_{11} - k_{22})^2-4(k_{11} d_{22} + k_{22} d_{11}) (d_{11}d_{22} - d_{12}^2) (d_{11} + d_{22})}{4(d_{11}+ d_{22})^2}.\nn
\ea

Suppose now that the damping force decreases in a uniform way, so we put
 $d_{11} = \varepsilon d_{11}'$, $d_{12} = \varepsilon d_{12}'$, $d_{22} = \varepsilon d_{22}'$, where $d_{11}$, $d_{12}$, $d_{22}$ are constants and $\varepsilon \rightarrow 0$. Then, for the inequality \rf{bb2} we get
\be{bb3}
\nu^2 < \nu_{cr}^2:= \frac{(k_{11}-k_{22})^2}{4}-\frac{(d_{11}'-d_{22}')^2(k_{11} - k_{22})^2}{4(d_{11}'+ d_{22}')^2}.
\ee
But if there is no damping, we have to make use of condition $B$, which gives
\be{bb4}
\nu^2 < {\nu_0}^2:= \frac{(k_{11}-k_{22})^2}{4}=\left( \frac{{\rm tr}{\bf K}}{2}\right)^2-\det{\bf K}.
\ee
Obviously
\be{bb5}
{\nu_0}^2-\nu_{cr}^2=\frac{(d_{11}'-d_{22}')^2(k_{11} - k_{22})^2}{4(d_{11}'+ d_{22}')^2}=\left[\frac{2{\rm tr}{\bf K}{\bf D}-{\rm tr}{\bf K}{\rm tr}{\bf D}}{2{\rm tr}{\bf D}}\right]^2 \ge 0,
\ee
where the expressions written in terms of the invariants of the matrices involved \cite{K07a} are valid also without the restrictions on the matrices $\bf D$ and $\bf K$ that were adopted in \cite{Bo63,B55}.
For the values of $\frac{2{\rm tr}{\bf K}{\bf D}-{\rm tr}{\bf K}{\rm tr}{\bf D}}{2{\rm tr}{\bf D}}$ which are small with respect to $\nu_0$ we can approximately write \cite{Ki04,K05}
\be{bb6}
\nu_{cr}\simeq\nu_0-\frac{1}{2\nu_0}\left[\frac{2{\rm tr}{\bf K}{\bf D}-{\rm tr}{\bf K}{\rm tr}{\bf D}}{2{\rm tr}{\bf D}}\right]^2.
\ee
If $\bf D$ depends on two parameters, say $\delta_{1}$ and $\delta_{2}$, then \rf{bb6} has a canonical form \rf{w} for the Whitney's umbrella in the $(\delta_1,\delta_2,\nu)$-space.
Due to discontinuity existing for  $2{\rm tr}{\bf K}{\bf D}-{\rm tr}{\bf K}{\rm tr}{\bf D}\ne0$ the equilibrium may be stable if there is no
damping, but unstable if there is damping, however small it may be. We observe also that the critical non-conservative parameter, $\nu_{cr}$, depends on the \textit{ratio} of the damping coefficients and thus is strongly sensitive to the distribution of damping similarly to how it happens in other applications, including the viscous Chandrasekhar-Friedman-Schutz (CFS) instability of the Maclaurin spheroids \cite{LD1977}.

The analytical approximations of the type \rf{bb6} for the onset of the flutter instability in the general finite-dimensional and infinite-dimensional cases were obtained for the first time in the works \cite{K03a,K03b, Ki04, K05, KS05b, K07b} as a result of further development of the sensitivity analysis of simple and multiple eigenvalues in multiparameter families of non-self-adjoint operators. The previous important works include \cite{AY75,BBM89a,Br94,C95,Ha92,DS1979,J88,MK91,OR96,S1980-2,S1980-3,S96,SKM05,Sey}. Recent results on the perturbation analysis of dissipation-induced instabilities and the destabilization paradox are summarized in the works \cite{K2013dg} and \cite{LFA2016}.

 \section{An umbrella without dynamics}  \label{Wh}
 Part of global analysis, a topic of pure mathematics, is concerned with singularity theory, which deals with the geometric
 characterisation and classification of singularities (stationary points) of vector fields. In dynamics these singularities
 are recognised as equilibria of dynamical systems. Well-known representatives of this singularity school are
 Ren\'e Thom \cite{Th}
 and Christopher Zeeman. Among pure mathematicians they were exceptional as they promoted singularity theory
 as useful for real-life problems in biology, the social sciences and physics. Unfortunately their approach gave
 singularity theory a bad name as in their examples they used geometric methods without explaining a possible relation
 between realistic vector fields and dynamics. It makes little
 sense to describe equilibria and transitions (bifurcations) between equilibria without discussing actual causes that are tied in
 with dynamical processes and corresponding equations of motion. We want to stress here that
 notwithstanding the lack of dynamics the geometry of
 singularities as an ingredient of dynamical systems theory can be very useful.

\begin{figure}[h]
\begin{center}
\includegraphics[width=0.5\textwidth]{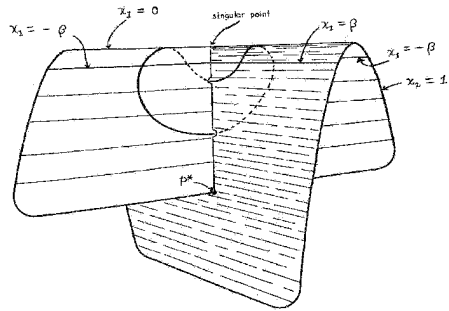}
\end{center}
\caption{\label{umb} Whitney's  umbrella, lowest dimensional case with 3 parameters \cite{W43,W44}.}
\end{figure}

 Before Ziegler's results a geometric result in singularity theory was obtained (1943-44) by Hassler Whitney.
 This result turned out to be
 an  excellent complement to Bottema's analytic approach.
 In his paper  \cite{W43}, Whitney described singularities of maps of a differential
$n$-manifold into $E^m$ with $m= 2n-1$. It turns out that in this case a special kind of
singularity plays a prominent role. Later, the local geometric structure of the manifold
near the singularity has been aptly called `Whitney's umbrella'. In Figure~\ref{umb} we reproduce a
sketch of the singular surface.

The paper \cite{W43} contains two main theorems. Consider the $C^k$ map $f: E^n \mapsto E^m$ with
$m= 2n-1$.
\begin{enumerate}
\item The map $f$ can be altered slightly, forming $f^*$, for which the singular points are
isolated. For each such an isolated singular point $p$, a technical regularity condition $C$ is valid
which relates to the map $f^*$ of the independent vectors near $p$ and of the differentials, the vectors
in tangent space.
\item Consider the map $f^*$ which satisfies condition $C$. Then we can choose coordinates
$x= (x_1, x_2, \cdots, x_n)$ in a neighborhood of $p$ and coordinates $y= (y_1, y_2, \cdots, y_m) $
(with $m= 2n-1$) in a neighborhood of $y=f(p)$ such that in a neighborhood of $f^*(p)$ we have
exactly
\begin{eqnarray*}
y_1 &=& x_1^2,\\
y_i &=& x_i,\,\,\,i=2, \cdots, n,\\
y_{n+i-1} &=& x_1x_i,\,\,\,i=2, \cdots, n.
\end{eqnarray*}
\end{enumerate}
If for instance $n=2,\, m=3$, the simplest interesting case, we have near the origin
\begin{equation}
y_1 = x_1^2,\,\,y_2 = x_2,\,\,y_3 = x_1x_2,
\end{equation}
so that $y_1 \geq 0$ and on eliminating $x_1$ and $x_2$:
\be{w} y_1y_2^2 - y_3^2 =0. \ee
Starting on the $y_2$-axis for $y_1=y_3=0$, the surface widens up for increasing values of $y_1$.
For each $y_2$, the cross-section is a parabola;
as $y_2$ passes through $0$, the parabola degenerates to a half-ray, and opens
out again (with sense reversed); see Figure~\ref{umb}.

Note that because of the $C^k$ assumption for the differentiable map $f$, the analysis is delicate.
There is a considerable simplification of the treatment if the map is analytical.

The analysis of singularities of functions and maps is a fundamental ingredient for bifurcation studies
of differential equations. After the pioneering work of Hassler Whitney and Marston Morse, it has become
a huge research field, both in theoretical investigations and in applications. We can not even present
a summary of this field here, so we restrict ourselves to citing a number of survey texts and
discussing a few key concepts and examples. In particular we mention \cite{AA88,Ar71,Ar83,AA93,GS85,GS88,G1997,F1992}. A monograph relating bifurcation theory with normal forms and
numerics is \cite{Kuz04}.

\section{Hopf bifurcation near 1:1 resonance and structural stability}

A study of the stability of equilibria of dynamical systems will usually involve the analysis of matrices obtained by
 linearisation of the
 equations of motion in a neighbourhood of the equilibria. This triggered off the study of {\em structural stability of matrices}
 as an independent topic in singularity theory \cite{Ar71, Ar83}.

More explicitly, consider a dynamical system described by the autonomous ODE
\[ \dot{\bf x} = {\bf f}({\bf x},{\bf p}),\,\,{\bf x} \in {\mathbb{R}}^n,\,\,{\bf f}:  {\mathbb{R}}^n \mapsto {\mathbb{R}}^n,\]
where ${\bf p}\in {\mathbb{R}}^k$ is a vector of parameters.
An equilibrium ${\bf x}_0$ of the system arises if ${\bf f}({\bf x}_0,{\bf p})={\bf 0}$. With a little smoothness of the
map ${\bf f}$ we can linearise near ${\bf x}_0$ so that we can write
\be{dyns}
\dot{\bf x} = {\bf A}_{\bf p}({\bf x}-{\bf x}_0) + {\bf g}({\bf x},{\bf p})
\ee
with ${\bf A}_{\bf p}$ a constant $n \times n -$ matrix, ${\bf g}({\bf x},{\bf p})$ contains higher-order terms only. In other words
\[ \lim_{{\bf x} \rightarrow {\bf x}_0}  \frac{\| {\bf g}({\bf x},{\bf p}) \|}{\| {\bf x} - {\bf x}_0 \|} = 0, \]
${\bf g}({\bf x},{\bf p})$ is tangent to the linear map in ${\bf x}_0$.
The properties of the matrix ${\bf A}_{\bf p}$ determine in a large number of cases the local behavior of the dynamical system.

Suppose that for ${\bf p} ={\bf 0}$, ${\bf A}_{\bf 0}$ has two equal non-zero imaginary eigenvalues and their complex conjugates, $\pm i \omega, \omega > 0$, and no other eigenvalues with zero real part. This equilibrium is called
a $1:1$ resonant double Hopf point \cite{G1997}. (Similarly, in a Hopf point the matrix of linearization has a single conjugate pair of imaginary eigenvalues $\pm i \omega$ and in a double Hopf point there are two distinct such pairs: $\pm i\omega_1$, $\pm i\omega_2$ \cite{G1997}.)
Then, without loss of generality, we may assume that the system \rf{dyns} has been already reduced to a centre manifold of dimension $n=4$. Considering further a generic case of double non-semisimple eigenvalues with geometric multiplicity 1, after a linear change of coordinates and re-scaling time to get $\omega=1$, we can transform ${\bf A}_{\bf 0}$ to \cite{F1992}
\be{tm}
\left(
  \begin{array}{rrrr}
    0 & -1 & 1 & 0 \\
    1 & 0 & 0 & 1 \\
    0 & 0 & 0 & -1 \\
    0 & 0 & 1 & 0 \\
  \end{array}
\right).
\ee
Setting $z_1=\Delta x_1+i\Delta x_2$ and $z_2=\Delta x_3+i\Delta x_4$, where $i=\sqrt{-1}$ and $\Delta {\bf x}={\bf x}-{\bf x}_0$, and assuming ${\bf A}_{\bf 0}$ to be in the form \rf{tm} we re-write \rf{dyns} at ${\bf p}={\bf 0}$ in the complex form \cite{F1992}
\be{cf}
\left(
  \begin{array}{c}
    \dot z_1 \\
    \dot z_2 \\
  \end{array}
\right)=\left(
          \begin{array}{cc}
            i & 1 \\
            0 & i \\
          \end{array}
        \right)\left(
  \begin{array}{c}
     z_1 \\
     z_2 \\
  \end{array}
\right)+\widetilde{{\bf g}}(z_1,z_2,\overline{z}_1,\overline{z}_2).
\ee
The second pair of equations governing the conjugates $\overline{z}_1$, $\overline{z}_2$ is omitted here for simplicity.

Arnold \cite{Ar71,Ar83} has proven that a universal unfolding of the linear vector field with the matrix $$\left(
                                                                                                       \begin{array}{cc}
                                                                                                         i & 1 \\
                                                                                                         0 & i \\
                                                                                                       \end{array}
                                                                                                     \right)
$$
is given by the three-parameter family of complex matrices
\be{unf}
\left(
  \begin{array}{cc}
    i+\alpha & 1 \\
    \mu_1+i \mu_2 & i +\alpha \\
  \end{array}
\right),
\ee
where $\alpha$, $\mu_1$, and $\mu_2$ are real parameters and versality is understood with respect to the group of similarity transformations
and a real positive scaling. The set of matrices with a resonant Hopf pair is a group orbit \cite{G1997}.

\begin{figure}
\begin{center}
\includegraphics[width=0.8\textwidth]{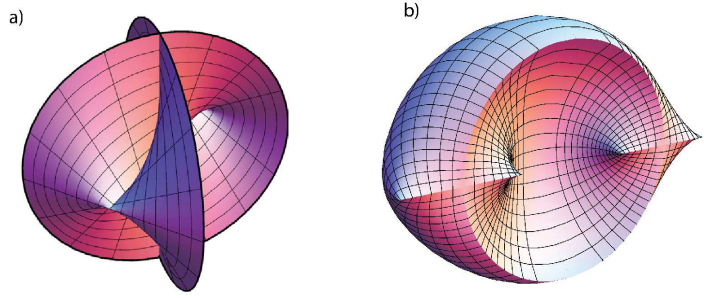}
\end{center}
\caption{\label{hk} (a) The Pl\"ucker conoid in the unfolding of a semisimple $1:1$ resonance has a pair of Whitney umbrellas \cite{HK10}; (b) two pairs of Whitney umbrellas on the boundary of the stability domain of a general 4-degrees-of-freedom dynamical system near a semisimple $1:1$ resonance  \cite{HK14}.}
\end{figure}

The universal unfolding has a pure imaginary eigenvalue if and only if there exists a real
number $\delta$ such that $(i(\delta - 1) - \alpha)^2 - (\mu_1 + i\mu_2) = 0$. Eliminating $\delta$ yields \cite{G1997}
$$\alpha^2(−4\mu_1 + 4\alpha^2) = \mu_2^2.$$

Setting $y_2 = \alpha$, $y_3 = \mu_2$, and $y_1 = −4\mu_1 + 4\alpha^2$ we reduce the equation to the form $y_3^2 = y_1y_2^2$, which is nothing else but
the normal form \rf{w} for the Whitney umbrella. The double Hopf points of \rf{unf} form the half-line $\alpha = \mu_2 = 0$, $\mu_1 = 0$. Along the
continuation $\mu_1 > 0$ of this half-line the eigenvalues of \rf{unf} are given by $\lambda = i\pm \sqrt{\mu_1}$. We see that the double Hopf points have codimension 2 and the resonant double Hopf points are of codimension 3.

If a family of matrices ${\bf A}({\bf p}) ={\bf A}(p_1, p_2, p_3, \ldots, p_k)$ has a $1:1$ resonant
double Hopf point, the universality of the unfolding \rf{unf} means that there exist smooth functions
$\alpha({\bf p})$, $\mu_1({\bf p})$, $\mu_1({\bf p})$, such that the Hopf structure of ${\bf A}({\bf p})$ near the $1:1$ resonant point
is the same as the Hopf structure of the unfolding with $\alpha$, $\mu_1$, $\mu_2$ replaced by $\alpha({\bf p})$, $\mu_1({\bf p})$, $\mu_2({\bf p})$.

Therefore, the stratified set of Hopf points in the neighborhood of a \textit{non-semisimple} $1:1$
resonance is a Whitney umbrella in $\bf p$-space too \cite{G1997}. The functions $\alpha({\bf p})$, $\mu_1({\bf p})$, $\mu_2({\bf p})$ can be found
approximately as truncated Taylor series with respect to the components of the vector $\bf p$ of the parameters \cite{Sey}.

Similar stratification of Hopf points near a \textit{semisimple} $1:1$ resonance involves pairs of Whitney umbrellas forming a Pl\"ucker conoid \cite{HK10,HK14}, see Figure~\ref{hk}.

Hoveijn and Ruijgrok (1995) were the first who applied these ideas to a practical problem exhibiting the Ziegler paradox. Namely, they considered a problem of widening due to dissipation of the zones of
the combination resonance \cite{YS75} in a system of two parametrically forced coupled oscillators \cite{HR95}.
The system models a rotating disk with oscillating suspension point introduced in \cite{RTV93}. Its linearized equations are
\ba{hrsys}
\ddot{x} + 2 \Omega \dot{y} + (1+ \varepsilon \cos \omega_0t)x + 2 \varepsilon \mu \dot{x} &=&0,\nn \\
\ddot{y} + 2 \Omega \dot{x} + (1+ \varepsilon \cos \omega_0t)y + 2 \varepsilon \mu \dot{y} &=&0.
\ea
It is assumed that for $\varepsilon=0$ the system \rf{hrsys} has two pairs of imaginary eigenvalues $\pm i \omega_1$, $\pm i \omega_2$ that depend on the parameter $\Omega$ representing the speed of rotation.  Of special interest is the case of the sum resonance
$\omega_0 = \omega_1+ \omega_2$.

Let parameters $\delta_1$ and $\delta_2$ control the detuning of the frequencies $\omega_1$ and $\omega_2$; then $\delta_+=\delta_1+\delta_2$ and $\delta_-=\delta_1-\delta_2$. The parameter $\delta_+$ is small and represents the detuning of the exact sum resonance: $\omega_0 = \omega_1+ \omega_2+\delta_+$ where $\delta_+=0$. Parameters $\mu_1$ and $\mu_2$ control the detuning of the damping from $\mu$; $\mu_+=\mu_1+\mu_2$, $\mu_-=\mu_1-\mu_2$.

The original nonlinear system that has the linearization \rf{hrsys} at zero detuning can be reduced to the following type of equation \cite{HR95}
\begin{equation} \label{HReq}
\dot{\bf z} = {\bf A}_0 {\bf z} + \varepsilon {\bf f}({\bf z}, \omega_0t; {\bf p}), \quad {\bf z}\in {\mathbb{R}}^4,
\end{equation}
where ${\bf A}_0$ is a $4 \times 4$ matrix with the eigenvalues $\pm i \omega_1, \pm i \omega_2$.
The vector of parameters ${\bf p}=(\delta_+,\delta_-,\mu_+,\mu_-)$ is used to control detuning from resonance and damping.

The
vector-valued function $\bf f$ is $2\pi$-periodic in $\omega_0 t$ and ${\bf f}(0, \omega_0 t; {\bf p}) = {\bf 0}$ for all $t$
and $\bf p$.
Since the origin is a stationary point of \rf{HReq}, one may ask how its stability
depends on the parameters. Analogous to the case of a single forced
oscillator, one can make a planar stability diagram by varying the strength $\varepsilon$
and the frequency $\omega_0$ of the forcing while fixing the other parameters. Also
in this case one obtains a resonance tongue at $\omega_0=\omega_1+\omega_2$.
However if damping is varied, the planar stability diagram does not change
continuously \cite{RTV93}.
For instance, applying zero damping ($\mu =0$, no damping detuning)
we find instability of the trivial solution (equilibrium) if $$|\delta_+| \leq 1.$$
\noindent
For $\mu >0$ the trivial solution is unstable if
\[ |\delta_+| \leq \omega_0 \sqrt{ \frac{1}{4}- \frac{\mu^2}{\omega_0^2}}. \]
\noindent
The instability interval depends discontinuously on damping coefficient $\mu$!

\begin{figure}
\begin{center}
\includegraphics[width=0.5\textwidth]{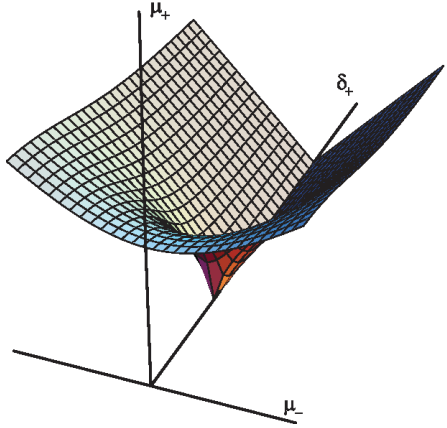}
\end{center}
\caption{\label{hrw} The critical surface for the damped combination resonance in $(\mu_+,\mu_-,\delta_+)$ space, where $\mu_+=\mu_1+\mu_2$, $\mu_-=\mu_1-\mu_2$, $\delta_+=\delta_1+\delta_2$. The parameters $\delta_1$ and $\delta_2$ control the detuning of the frequencies $\omega_1$ and $\omega_2$, the parameters $\mu_1$ and $\mu_2$ the damping of the oscillators. The self-intersection of the surface with the Whitney umbrella singularity is along the $\delta_+$ axis \cite{HR95}. }
\end{figure}

Hoveijn and Ruijgrok \cite{ HR95} presented a geometrical explanation of this dissipation-induced instability using `all' the parameters as unfolding parameters first putting the equation \rf{HReq} into a normal form \cite{Ar83,SVM} which is similar to that of the non-semisimple $1:1$ resonance studied in \cite{F1992}.

In the normalized equation the time dependence appears only in
the high order terms. But the autonomous part of this equation contains
enough information to determine the stability regions of the origin.

The second step was to test the linear autonomous part ${\bf A}({\bf p})$ of the normalized
equation for structural stability. This family
of matrices is parameterized by the detunings of the frequencies $\omega_1$ and $\omega_2$ and of the
damping parameter $\mu$.

Identifying the most degenerate member of
this family one can show that ${\bf A}({\bf p})$ is its versal unfolding in the sense of
Arnold \cite{Ar71,Ar83}.
Put differently, the family ${\bf A}(\delta_+,\delta_-,\mu_+,\mu_-)$ is structurally stable, whereas
${\bf A}(\delta_+,\delta_-,0,0)$ is not. Therefore the stability diagram actually `lives' in a
four dimensional space. In this space, the stability regions of the origin are
separated by a critical surface which is the hypersurface where ${\bf A}({\bf p})$ has at
least one pair of imaginary complex conjugate eigenvalues. This
critical surface is diffeomorphic to the Whitney umbrella, see Figure~\ref{hrw}.

It is the
singularity of the Whitney umbrella that causes the discontinuous behaviour
of the planar stability diagram for the combination resonance in the presence of dissipation. The structural stability argument guarantees
that the results are `universally valid' and qualitatively hold for
every system in sum resonance. For technical details we refer again to \cite{HR95}.

 \section{Abscissa minimization, robust stability and heavy damping}

\begin{figure}
\begin{center}
\includegraphics[width=0.99\textwidth]{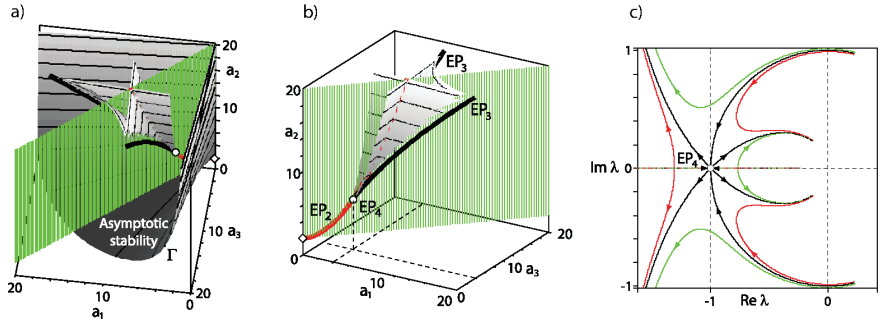}
\end{center}
\caption{\label{ws} (a) A singular boundary $\Gamma$ of the domain of asymptotic
stability \rf{domain} of the polynomial \rf{polws} with the Whitney umbrella
singularity at the point $(a_1, a_3, a_2) = (0, 0, 2)$, marked by the
diamond symbol. (b) The tangent cone, the EP-sets, and the discriminant surface with the Swallowtail singularity at $EP_4$. The domain of heavy damping is inside the `spire'. (c) Trajectories of roots of the polynomial \rf{polws} when $a_2$
increases from 0 to 15 and: $a_1 = a_3 = 4$ (black); $a_1 = 4$, $a_3 = 3.9$ (red);
$a_1 = 3.9$, $a_3 = 4$ (green). The global minimum of the abscissa is attained
when all the roots coalesce into the quadruple root $\lambda = -1$ ($EP_4$) \cite{KO2013}.}
\end{figure}

Let us return to the work of Bottema \cite{B56}. The conditions
\be{domain}
a_1 > 0,\quad a_3 > 0, \quad a_2 > 2 + \frac{(a_1 - a_3)^2}{a_1a_3}
>0
\ee
are necessary and sufficient for the polynomial
\be{polws}
p(\lambda) = \lambda^4 + a_1\lambda^3 + a_2\lambda^2 + a_3\lambda +1
\ee
to be Hurwitz. The domain \rf{domain} was plotted by Bottema in the
$(a_1, a_3, a_2)$-space, Figure~\ref{ws}a.

A part of the plane $a_1=a_3$ that lies inside the domain of asymptotic stability constitutes a set of all directions leading from the point (0, 0, 2) to the stability region
\be{tc}
\{(a_1, a_3, a_2):~~ a_1 = a_3,~~ a_1 > 0,~~ a_2 > 2\}.
\ee
The \textit{tangent cone} \rf{tc} to the domain of asymptotic stability at
the Whitney umbrella singularity, which is shown in green in
Figure~\ref{ws}a,b, is degenerate in the $(a_1, a_3, a_2)$ -- space because it selects a
set of measure zero on a unit sphere with the center at the singular
point \cite{L80,L82}.

The singular point $(a_1, a_3, a_2) = (0, 0, 2)$ corresponds to a
double complex-conjugate pair of roots $\lambda = \pm i$ of the polynomial \rf{polws}.
The fact that multiple roots of a polynomial
are sensitive to perturbation of the coefficients is a phenomenon that was studied already by Isaac Newton,
who introduced the so-called Newton polygon to determine the
leading terms of the perturbed roots as fractional powers of
a perturbation parameter. It follows that, in matrix analysis, eigenvalues are in general not locally Lipschitz at points
in matrix space with non-semi-simple eigenvalues \cite{Bu06},
and, in the context of dissipatively perturbed Hamiltonian
systems, \cite{MO95}. Thus, it has been well-understood for
a long time that perturbation of multiple roots or multiple
eigenvalues on or near the stability boundary is likely to lead to
instability \cite{KS05a}.

Because of the sensitivity of multiple roots and eigenvalues to
perturbation, in engineering and control-theoretical applications
a natural desire is to ``cut the singularities off'' by constructing
convex inner approximations to the domain of asymptotic
stability. Nevertheless, multiple roots per se are not undesirable.

Indeed, multiple roots also occur deep inside the domain
of asymptotic stability. Although it might seem paradoxical at first
sight, such configurations are actually obtained by minimizing the
so-called \textit{polynomial abscissa} in an effort to make the asymptotic
stability of a linear system more \textit{robust}, as we now explain.
\\\\
\noindent
{\bf Abscissa minimization and multiple roots}\\
The abscissa of a polynomial $p(\lambda)$ is the maximal real part of its
roots \cite{BGMO2012}:
\be{ad}
a(p) = \max\{{\rm Re}~ \lambda : ~~p(\lambda) = 0\}.
\ee
We restrict our attention to monic polynomials with real coefficients and fixed degree $n$: since these have $n$
free coefficients, this space is isomorphic to $\mathbb{R}^n$. On this space, the abscissa is a continuous but non-smooth, in fact non-Lipschitz, as well as non-convex, function whose variational properties have been extensively studied using non-smooth variational analysis \cite{BGMO2012,Bu06}.

Now set $n = 4$, consider the set of polynomials $p(\lambda)$ defined in
\rf{polws}, and consider the restricted set of coefficients
\be{setco}
\left\{(a_1, a_3, a_2) : ~~a_1 = a_3,~~ a_2 = 2 + \frac{a_1^2}{4}\right\}
\ee
On this set the roots are
$$
\lambda_1 = \lambda_2 = -\frac{a_1}{4}-\frac{1}{4}\sqrt{a_1^2-16},\quad
\lambda_3 = \lambda_4 = -\frac{a_1}{4}+\frac{1}{4}\sqrt{a_1^2-16}.
$$
When $0 \le a_1 < 4$ $(a_1 > 4)$, the roots $\lambda_{1,2}$ and $\lambda_{3,4}$ are complex
(real) with each pair being double, that is with multiplicity two. At
$a_1 = 4$ there is a quadruple real eigenvalue $-1$. So, we refer to the
set \rf{setco} as a set of exceptional points \cite{K2013dg,KI2017,KI2018} (abbreviated as the EP-set).

When $a_1 > 0$, the EP-set \rf{setco} (shown by the red curve in
Figure~\ref{ws}a,b) lies within the tangent cone \rf{tc} to the domain of
asymptotic stability at the Whitney umbrella singularity $(0, 0, 2)$.
The points in the EP-set all define polynomials with two double
roots (denoted $EP_2$) except $(a_1, a_3, a_2) = (4, 4, 6)$, at which $p$ has
a quadruple root and is denoted $EP_4$; see Figure~\ref{ws}a,b.

Let us consider how the roots move in the complex plane when $a_1$ and $a_3$ coincide and are set
to specific values and $a_2$ increases from zero, as shown by black curves in Figure~\ref{ws}c.
When $a_1 = a_3 < 4$, the roots that initially have positive real parts and thus
correspond to unstable solutions move along the unit circle to the left, cross the imaginary axis at
$a_2 = 2$ and merge with another complex conjugate pair of roots at $a_2 = 2 + \frac{a_1^2}{4}$,
i.e., at the EP-set. Further increase in $a_2$ leads to the splitting of
the double eigenvalues, with one conjugate pair of roots moving back toward the imaginary axis.
By also considering the case $a_1 = a_3 > 4$, it is clear that when $a_1$ and $a_3$ coincide, the choice
$a_2 = 2 + \frac{a_1^2}{4}$ minimizes the abscissa, with the polynomial $p$ on
the EP-set.

Furthermore, when $a_1 = a_3$ is increased toward $4$ from below, the coalescence points ($EP_2$) move around
the unit circle to the left. This conjugate pair of coalescence
points merges into the quadruple real root $\lambda = -1$ ($EP_4$) when
$a_1 = a_3 = 4$ and hence $a_2 = 6$, as is visible in Figure~\ref{ws}c. If $a_1 = a_3$
is increased above 4 the quadruple point $EP_4$ splits again into
two exceptional points $EP_2$, one of them inside the unit circle.

Thus, all indications are that the abscissa is minimized by the
parameters corresponding to $EP_4$, with a quadruple root at $-1$.

In fact, application of the following theorem shows that the abscissa
of \rf{polws} is globally minimized by the $EP_4$ parameters.

\begin{theorem}(\cite{BGMO2012}, Theorems 7 and 14)

Let $\mathbb{F}$ denote either the real field $\mathbb{R}$ or the complex field $\mathbb{C}$.
Let $b_0$, $b_1$, $\ldots$, $b_n \in \mathbb{F}$ be given (with $b_1, \ldots, b_n$ not all zero)
and consider the following family of monic polynomials of degree $n$ subject to a single affine constraint on the
coefficients:
$$
P = \left\{\lambda^n + a_1\lambda^{n-1} + \ldots + a_{n-1}\lambda + a_n : \quad
b_0 + \sum_{j=1}^n b_ja_j = 0,\quad a_i \in \mathbb{F}\right\}.
$$
Define the optimization problem
\be{op15}
a^* := \inf_{p \in P} a(p).
\ee
Let
$$
h(\lambda) = b_n\lambda^n + b_{n-1}\left(
                                     \begin{array}{c}
                                       n \\
                                       n-1 \\
                                     \end{array}
                                   \right)
\lambda^{n-1} + \ldots + b_1\left(
                                     \begin{array}{c}
                                       n \\
                                       1 \\
                                     \end{array}
                                   \right)\lambda + b_0.
$$
First, suppose $\mathbb{F} = \mathbb{R}$. Then, the optimization problem \rf{op15} has the
infimal value
$$
a^*= -\max\left\{\zeta \in \mathbb{R} : ~~h^{(i)}(\zeta) = 0 ~~for~~some~~ i \in \{0, \ldots , k - 1\}\right \},
$$
where $h^{(i)}$ is the $i$-th derivative of $h$ and $k = \max\{j : ~~b_j \ne 0\}$.
Furthermore, the optimal value $a^*$ is attained by a minimizing polynomial $p^*$ if and only if $-a^* $ is a root of $h$ (as opposed to one of its
derivatives), and in this case we can take
$$
p^*(\lambda) = (\lambda - \gamma)^n \in P,\quad  \gamma = a^*.
$$
Second, suppose $\mathbb{F} = \mathbb{C}$. Then, the optimization problem \rf{op15}
always has an optimal solution of the form
$$
p^*(\lambda) = (\lambda - \gamma)^n \in P,\quad  {\rm Re}~ \gamma = a^*,
$$
with $-\gamma$ given by a root of $h$ (not its derivatives) with largest real part.
\end{theorem}

In our case, $\mathbb{F} = \mathbb{R}$, $n = 4$ and the affine constraint on the
coefficients of $p$ is simply $a_4 = 1$. So, the polynomial $h$ is given by
$h(\lambda) = \lambda^4 - 1$.

Its real root with largest real part is 1, and its derivatives have
only the zero root. So, the infimum of the abscissa $a$ over the
polynomials \rf{polws} is $-1$, and this is attained by
\be{minpo}
p^*(\lambda) = (\lambda + 1)^4 = \lambda^4 + 4\lambda^3 + 6\lambda^2 + 4\lambda + 1,
\ee
that is, with the coefficients at the exceptional point $EP_4$. There is
nothing special about $n = 4$ here; if we replace $4$ by $n$ we find that
the infimum is still $-1$ and is attained by
$$
p^*(\lambda) = (\lambda + 1)^n.
$$

\noindent
{\bf Swallowtail singularity as the global minimizer of the abscissa}

It is instructive to understand
the set in the $(a_1, a_3, a_2)$-space where the roots of the polynomial \rf{polws} are real and negative,
but not necessarily simple, which is given by the discriminant surface of the polynomial. A part of it is shown in
Figure~\ref{ws}a,b. At the point $EP_4$ with the coordinates $(4, 4, 6)$ in the $(a_1, a_3, a_2)$-
space the discriminant surface has the \textit{Swallowtail} singularity,
which is a generic singularity of bifurcation diagrams in three-parameter families of real matrices \cite{Ar71,Ar83}.

Therefore, the coefficients of the globally minimizing polynomial \rf{minpo} are exactly at the Swallowtail singularity of the
discriminant surface of the polynomial \rf{polws}.

In the region in side the ``spire'' formed by the discriminant
surface all the roots are simple real and negative. Owing to this
property, this region, belonging to the domain of asymptotic stability (see Figure~\ref{ws}a), plays an important role in stability
theory. Physical systems with semi-simple real and negative eigenvalues are called \textit{heavily damped}. The solutions of the heavily damped systems do not oscillate and monotonically decrease, which is favorable for applications in robotics and automatic control.

Now we can give the following interpretation of the Bottema stability
diagram shown in Figure~\ref{ws}a \cite{KO2013}. The dissipative system with the characteristic polynomial \rf{polws} is
asymptotically stable inside the domain \rf{domain}. The boundary of the domain \rf{polws} has the
Whitney umbrella singular point at $a_1 = 0$, $a_3 = 0$, and $a_2 = 2$.

The domain corresponding to heavily damped systems is confined between three hypersurfaces of the discriminant surface and has a form of a trihedral spire with the Swallowtail singularity at its cusp at $a_1 = 4$, $a_2 = 6$, and $a_3 = 4$.
The Whitney umbrella and the Swallowtail singular points are connected by the EP-set given by \rf{setco}. At the
Swallowtail singularity of the boundary of the domain of heavily damped systems, the abscissa of the characteristic polynomial
of the damped system attains its global minimum.

Therefore, by minimizing the spectral abscissa one finds points
at the boundary of the domain of heavily damped systems.
Furthermore, the sharpest singularity at this boundary corresponding to a quadruple real eigenvalue $\lambda =-1$ with
the Jordan block of order four is the very point where all the modes of the system with two degrees of freedom are decaying to zero as
rapidly as possible when $t \rightarrow \infty$.

\end{document}